\documentclass[12pt]{amsart}
\usepackage{amsmath,amsfonts,amssymb,amsthm,epsfig}
\usepackage{latexsym}
\usepackage{amssymb}

\usepackage{hyperref}

\usepackage{setspace}

\usepackage[margin=1.5in]{geometry}

\def\({\left(}
\def\){\right)}
\def\[{\left[}
\def\]{\right]}
\def\<{\left\langle}
\def\>{\right\rangle}

\newtheorem{theorem}{Theorem}[section]

\newtheorem{remark}[theorem]{Remark}

\DeclareMathOperator{\sign}{sign}

\begin{document}

\title[Plane partitions with 2-periodic weights]
{Plane partitions with 2-periodic weights}

\author[Sevak Mkrtchyan]{Sevak Mkrtchyan}
\email{sevakm@math.cmu.edu}
\address{Carnegie Mellon University, Pittsburgh, PA}

\begin{abstract}
We study scaling limits of skew plane partitions with periodic weights under several boundary conditions. We compute the correlation kernel of the limiting point process in the bulk and near turning points on the frozen boundary. The turning points that appear in the homogeneous case split in our model into pairs of turning points macroscopically separated by a ``semi-frozen'' region. As a result the point process at a turning point is not the GUE minor process, but rather a pair of GUE minor processes, non-trivially correlated. 

We also study an intermediate regime when the weights are periodic but all converge to $1$. In this regime the limit shape and correlations in the bulk are the same as in the case of homogeneous weights and periodicity is not visible in the bulk. However the process at turning points is still not the GUE minor process.
\end{abstract}

\maketitle

\tableofcontents

%

\section{Introduction}

\begin{figure}[ht]
\caption{\label{fig:boxes} An example of a skew plane partition with the corresponding stacks of cubes when $\lambda=\{3,2,2,1\}$.}
\includegraphics[width=8cm]{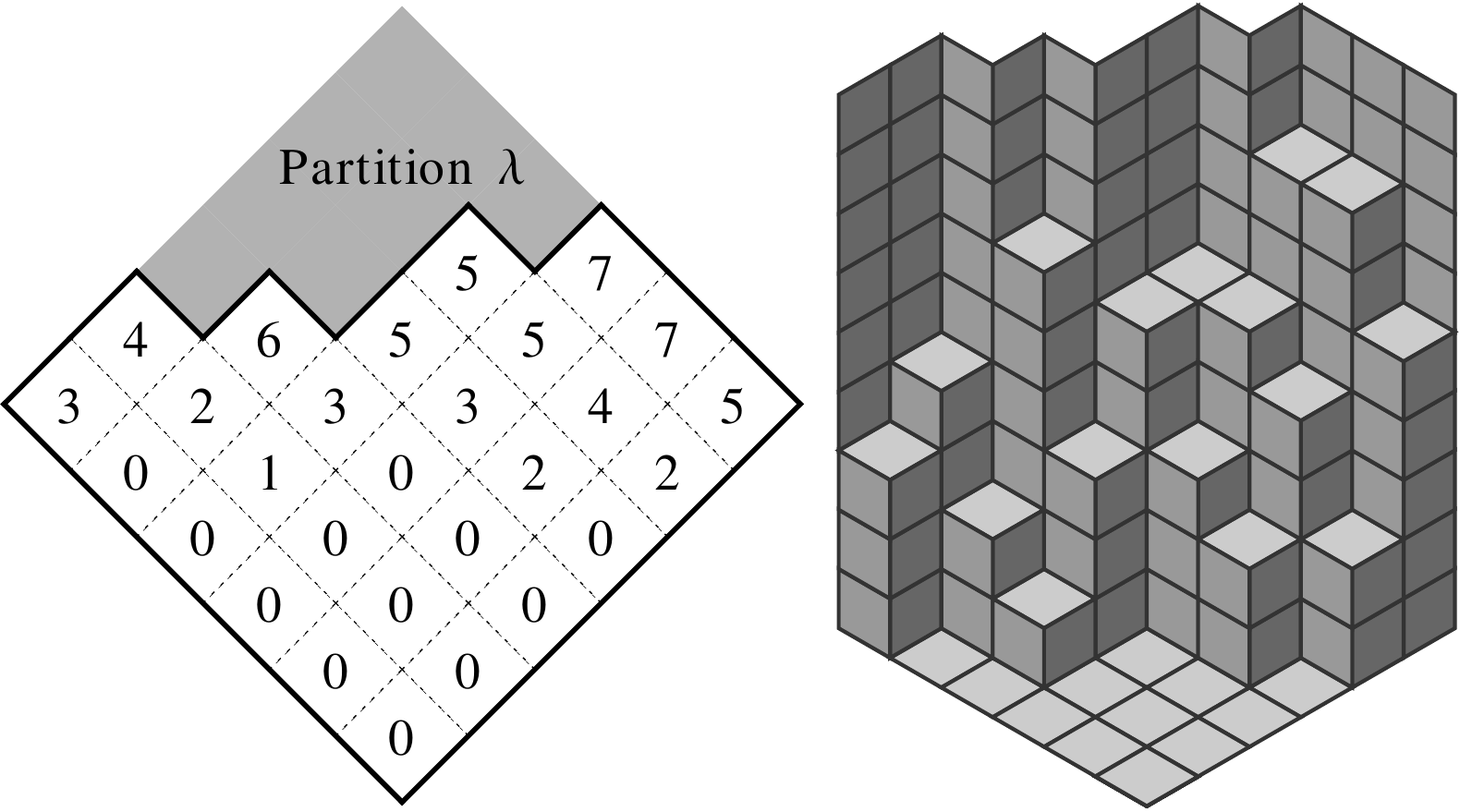}
\end{figure}

In this paper we study scaling limits of skew plane partitions with periodic weights under several boundary conditions. Recall that a \textit{partition} is a weakly decreasing sequence of non-negative integers where all but finitely many terms are zero. 
Given a partition $\lambda=(\lambda_1,\lambda_2,\dots)$, a \textit{skew plane partition} with boundary $\lambda$ confined to a $c\times d$ box is an array of non-negative integers $\pi=\{\pi_{i,j}\}$ defined for all $1\leq i\leq c$, $\lambda_i< j\leq d$, which are weakly decreasing in $i$ and $j$. We will denote the set of such skew plane partitions by $\Pi_\lambda^{c,d}$. 
We can visualize a skew plane partition $\pi$ as a collection of stacks of identical cubes where the number of cubes in position $(i,j)$ is equal to $\pi_{i,j}$, as shown in Figure \ref{fig:boxes}. 
The shape cut out by the partition $\lambda$ is often referred to as the ``back wall''.

Scaling limits of boxed plane partitions (plane partitions with restricted height and confined to a box) when all three dimensions grow at the same scale have been studied extensively in the mathematical literature. Limit shapes of random boxed plane partitions under the uniform measure were studied in e.g. \cite{CohnLarsenPropp},\cite{J}, \cite{KO}. The local correlations in the case of the uniform measure were studied in e.g. \cite{J},\cite{K1},\cite{G}. Boxed plane partitions under more general measures were studied in e.g. \cite{BorGorRainsqdistr}, \cite{BorGorShuffling}, \cite{BeteaElliptic}. Many ideas about limit shapes and the structure of correlations of the corresponding tiling models can be traced to the physics literature, e.g. \cite{NienHilBloTriangularSOS}.

A natural replacement for the uniform measure for skew plane partitions is the so-called ``volume'' measure: 
\begin{equation}
\label{eq:qvolume}
\mathbb{P}_q(\pi)\propto q^{|\pi|},
\end{equation}
where $q\in(0,1)$ and $|\pi|=\sum_{i,j}\pi_{i,j}$ is the volume of the plane partition $\pi$, i.e. the total number of stacked cubes. The limit shape phenomenon for such measures can be proven quite generally
\cite{CohnKenyonPropp}.

\subsubsection*{Schur process}
A skew plane partition can be represented as a sequence of interlacing Young diagrams consisting of its diagonal slices. More precisely, a skew plane partition $\pi\in\Pi^{c,d}_\lambda$ can be identified with the sequence of partitions $\{\pi(t)\}_{c< t< d}$ defined by 
\begin{equation*}
\pi(t)=(\pi_{i,t+i},\pi_{i+1,t+i+1},\dots),
\end{equation*}
where $i$ is such that $\lambda_j< t+j$ for $j\geq i$.
Under this identification the measure \eqref{eq:qvolume} corresponds to the measure
\begin{equation}
\label{eq:qvolumeSchurProc}
\mathbb{P}_q(\{(\pi(t)\}_{-c<t<d})\propto\prod_{-c<t<d}q^{|\pi(t)|}
\end{equation}
on sequences of partitions, were $|\pi(t)|$ is the size of the partition $\pi(t)$, i.e. the sum of all its entries. In this language it is natural to consider a non-homogeneous variant of the measure \eqref{eq:qvolumeSchurProc}: a measure where the value of the weight $q$ depends on the ``time'' coordinate $t$. More precisely, for a sequence $\bar{q}=\{q_t\}_{-c<t<d}$ of positive real numbers introduce a probability measure on the set $\Pi_\lambda^{c,d}$ by 
\begin{equation}
\label{eq:qinhomog}
\mathbb{P}_{\lambda,\bar{q}}(\pi)=\frac 1Z \prod_{-c<t<d}q_t^{|\pi(t)|},
\end{equation}
where $Z$ is the normalizing factor
\begin{equation*}
Z=\sum_{\pi\in\Pi_\lambda^{c,d}}\prod_{-c<t<d}q_t^{|\pi(t)|},
\end{equation*}
called the partition function. We refer to this as the inhomogeneous measure, or the measure with inhomogeneous weights, since in the interpretation of skew plane partitions as stacks of cubes it corresponds to giving different weights to cubes depending on which diagonal slice they fall on. 

Okounkov and Reshetikhin introduced a broad family of probability measures on sequences of partitions, called the Schur Process, of which the measure \eqref{eq:qinhomog} is a specialization \cite{OR1}. The Schur Process can be thought of as a time-dependent version of a measure on partitions introduced by Okounkov in \cite{O1}, called Schur measure. In our notation, the parameter $t$ plays the role of time.
In \cite{OR2} Okounkov and Reshetikhin showed that the  measures from the Schur Process are determinantal point processes, and gave a contour integral representation for the correlation kernel (see Theorem \ref{thm:fin-corr2} below). Using this they studied the scaling limit of random skew plane partitions with respect to the homogeneous measure \eqref{eq:qvolume} when the back wall approaches in the limit a piecewise linear curve with lattice slopes $\pm 1$ \cite{OR1},\cite{OR2}. Piecewise linear walls of non-lattice slopes were studied in \cite{BMRT}. Arbitrary piecewise linear back walls were considered in \cite{M}. 

\subsection{Main results}
While the correlation kernel of the underlying determinantal process in the case of inhomogeneous weights was computed in \cite{OR2}, the thermodynamic limit has never been studied. In this paper we carry such a study in the special case when the weights are periodic. For fixed $\alpha\geq1$ we define the weights $q_t$ by
\begin{equation}
\label{eq:periodicweights}
q_t=q_{t,r}=\left\{
\begin{array}{ll}
e^{-r}\alpha,&t\text{ is odd}\\
e^{-r}\alpha^{-1},&t\text{ is even}
\end{array}\right.,
\end{equation}
and study the scaling limit of the system when $r\rightarrow 0$. Note, that \eqref{eq:periodicweights} turns into the homogeneous measure \eqref{eq:qvolume} when $\alpha=1$.

In \cite{M} it was shown that for the homogeneous measure the scaling limit only depends on the macroscopic limit of the scaled back wall. In the case of the measure with weights \eqref{eq:periodicweights} when $\alpha\neq 1$, the system is very sensitive to microscopic changes to the back walls, and unlike the homogeneous case, the back walls have to be chosen carefully in order for the measure to be well defined (i.e. for the partition function to be finite). 
As a representative example, we study random skew plane partitions with the back wall given by a staircase, which in the scaling limit converges to a piecewise linear curve with three linear sections of slopes $1$, $0$ and $-1$ respectively.

\subsubsection{Bulk correlations and the frozen boundary} In Section \ref{sec:bulkCorrKer} we compute the asymptotics of the correlation kernel in the bulk. In the homogeneous case the distribution of the horizontal tiles in the neighbourhood of a point in the bulk converges
to a translation invariant ergodic Gibbs measure \cite{BMRT}, which were classified by Kenyon \cite{KenyonGibbs}. In the inhomogeneous case this is not the case any more, since when $\alpha>1$ the process is only $2\mathbb{Z}\times\mathbb{Z}$ translation invariant. 

In Section \ref{sec:frozenBoundary} We examine the frozen boundary in the three different regimes of interest:
\begin{itemize}
\item \textit{Unbounded floor:} The main difference between the inhomogeneous and homogeneous cases is that in the former the system is bounded everywhere except in the four tentacles that arise, whereas in the latter the system is always unbounded at linear sections of non-lattice slopes. In particular, the staircase shaped back wall at the linear section of slope $0$ stays frozen in the scaling limit in the inhomogeneous case.

\begin{figure}[ht]
\caption{\label{fig:unbddSample} The frozen boundary and an exact sample in the unbounded case. On the left side $\alpha=1$, while on the right side $\alpha>1$.}
\includegraphics[width=13.8cm]{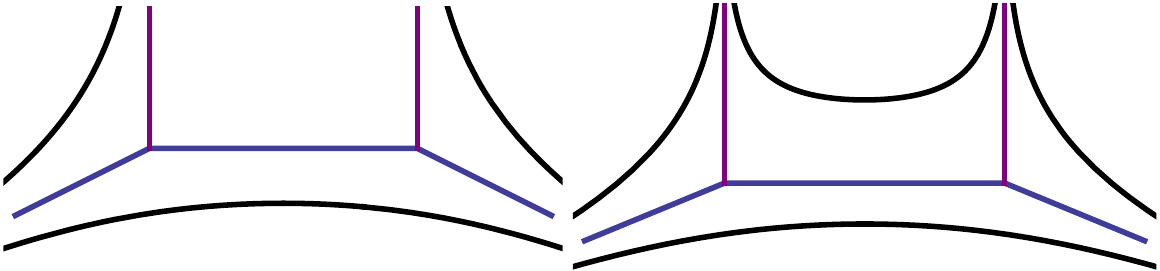}
\includegraphics[width=6.9cm]{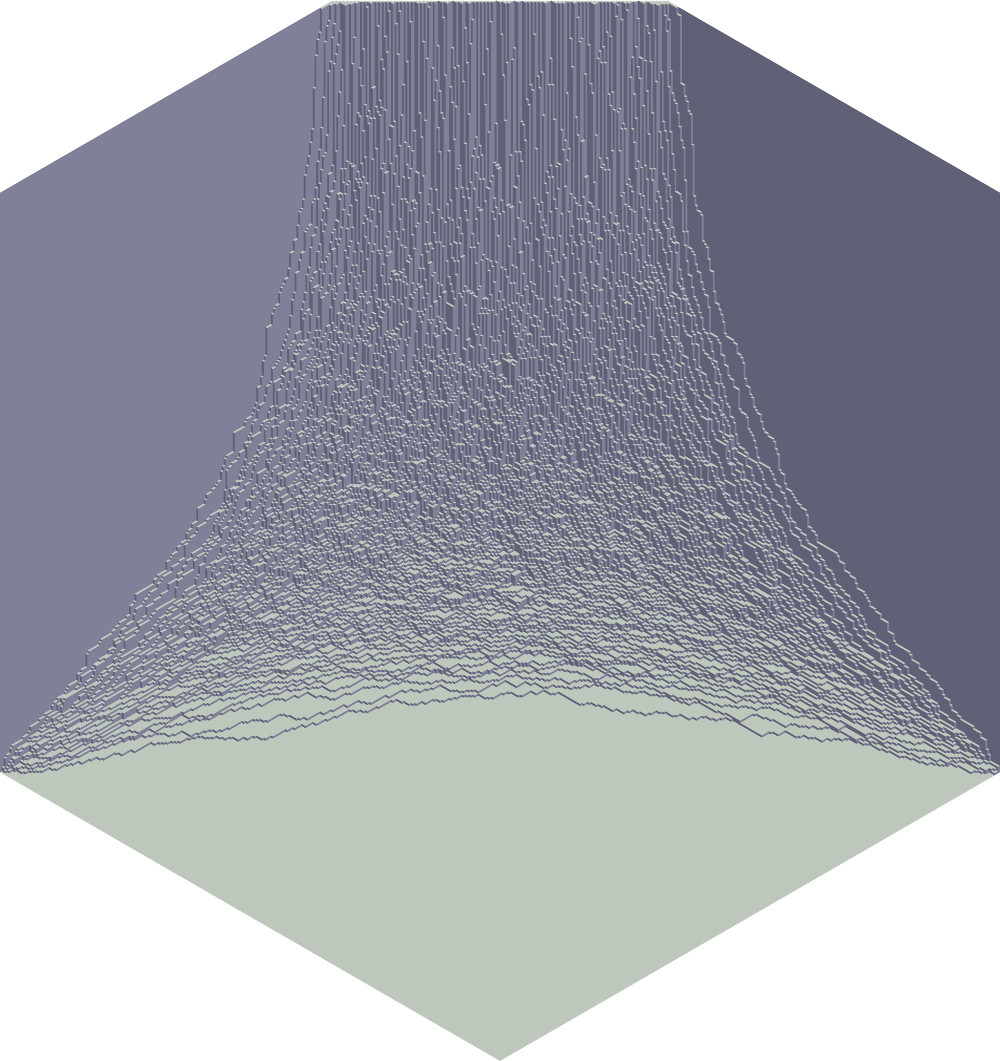}
\includegraphics[width=6.9cm]{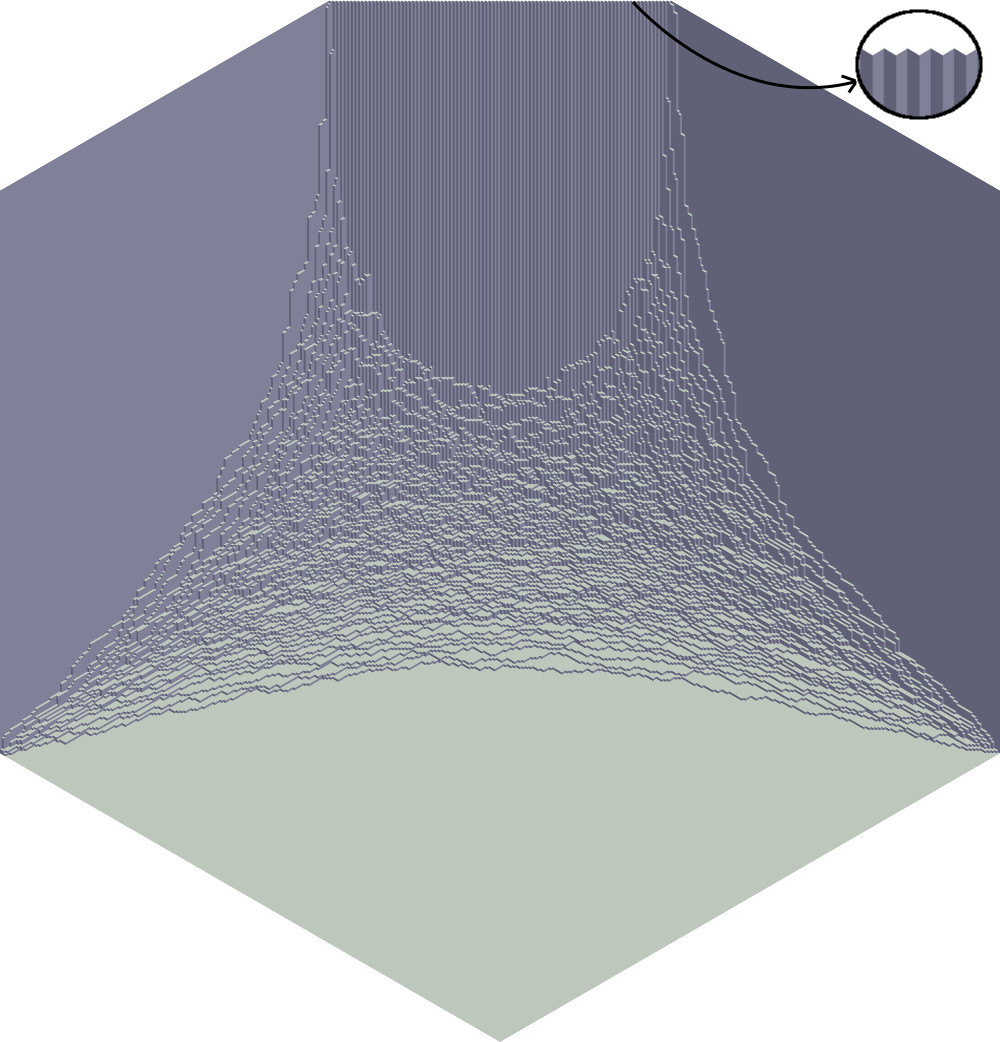}
\end{figure}

Note, that the shape of this frozen boundary can be predicted by methods from \cite{KOS}. A plot of the frozen boundary and an exact sample are shown in Figure \ref{fig:unbddSample}.
\item \textit{Bounded floor:} In this case, generically, the floor will be pentagonal in the scaling limit. The system develops four turning points near the vertical boundary, two on the left and two on the right. A notable difference from the homogeneous case is the nature of the turning points. We discuss this in detail in Section \ref{sec:turningPointsOverview}. A plot of the frozen boundary and an exact sample are shown in Figure \ref{fig:bddSample}.

\begin{figure}[ht]
\caption{\label{fig:bddSample} The frozen boundary and an exact sample in the bounded case. On the left side $\alpha=1$, while on the right side $\alpha>1$.}
\includegraphics[width=13.8cm]{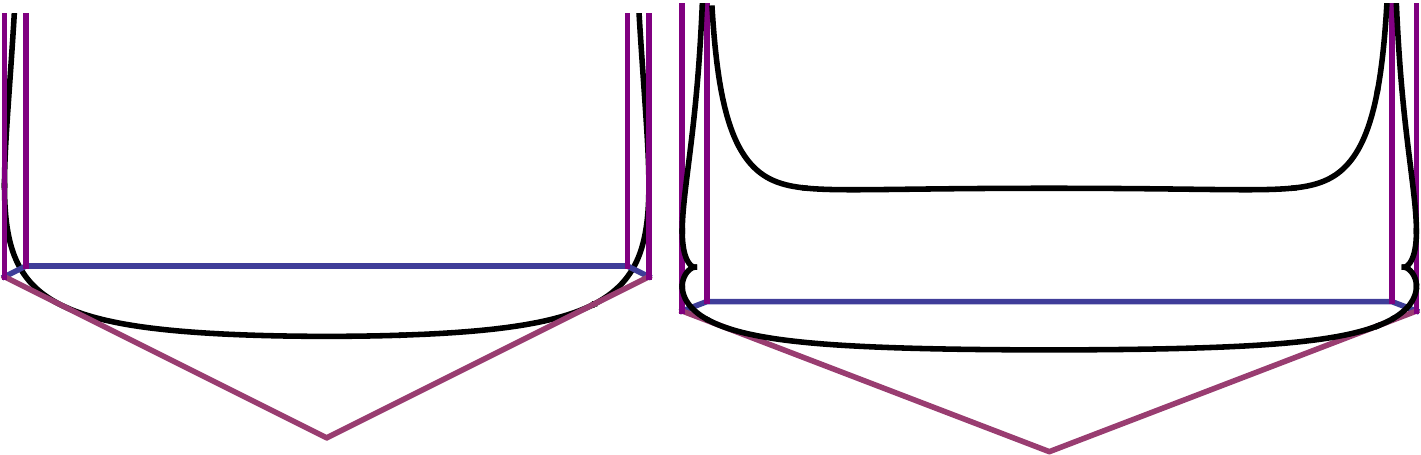}
\includegraphics[width=6.9cm]{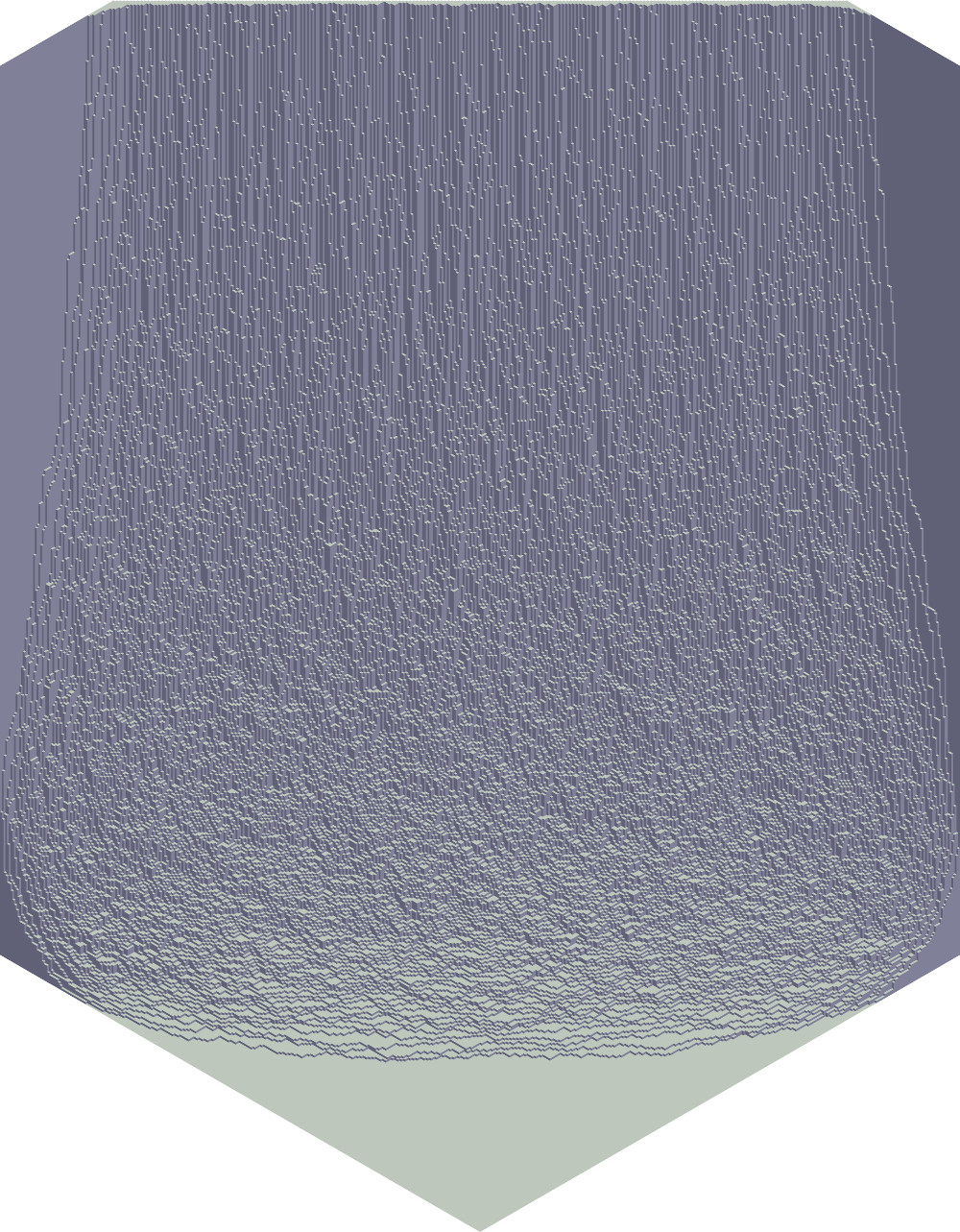}
\includegraphics[width=6.9cm]{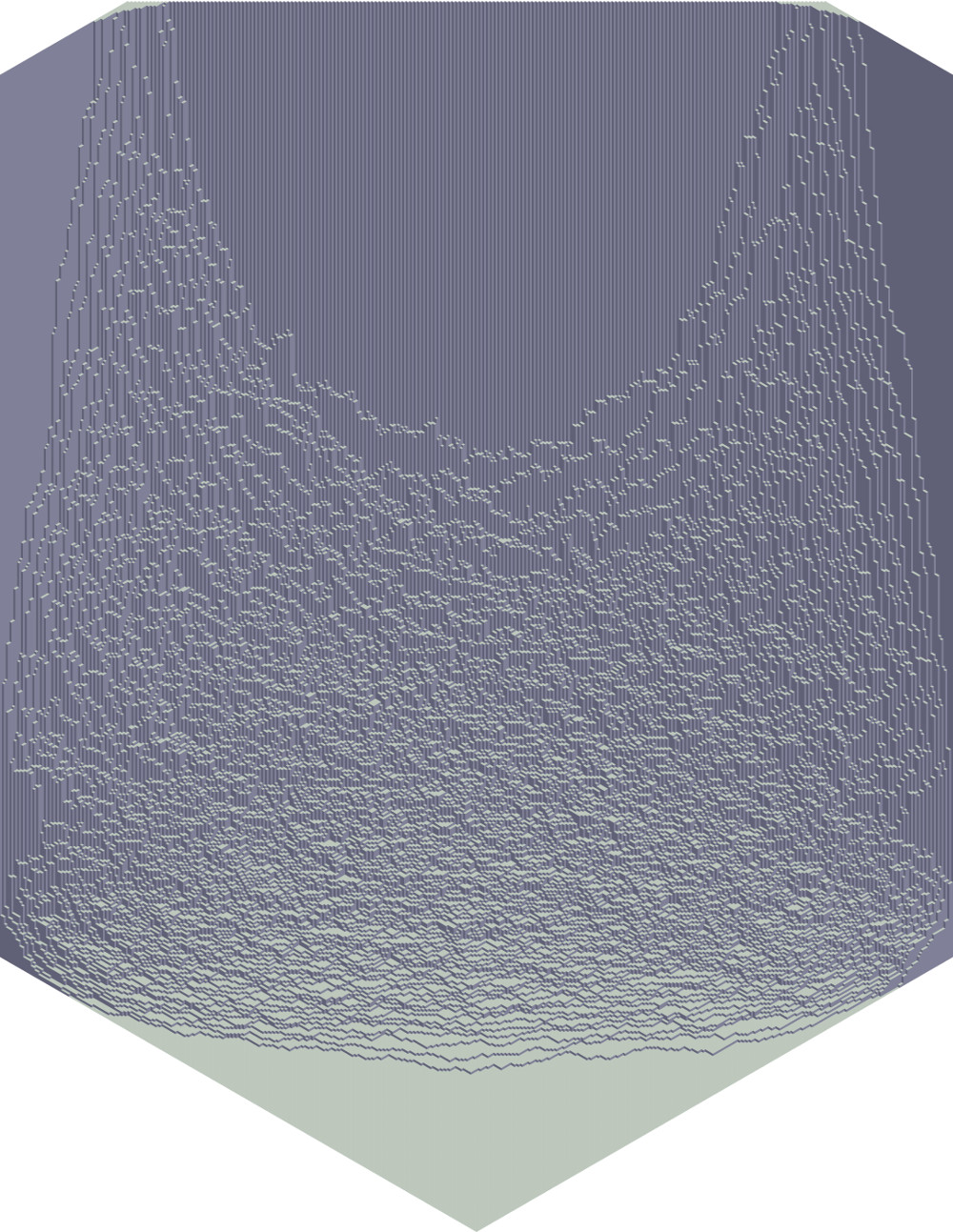}
\end{figure}

\item \textit{Triangular floor:} When the staircase back wall is as large as possible, the floor in the scaling limit is triangular. The corresponding limit for the homogeneous measure was studied in \cite{BMRT}, where it was shown that the disordered region is infinitely tall everywhere. In contrast, we find that when $\alpha>1$, the disordered region is bounded. A plot of the frozen boundary and an exact sample are shown in Figure \ref{fig:triangleSample}.

\begin{figure}[ht]
\caption{\label{fig:triangleSample} The frozen boundary and an exact sample in the triangular case. On the left side $\alpha=1$, while on the right side $\alpha>1$.}
\includegraphics[width=13.8cm]{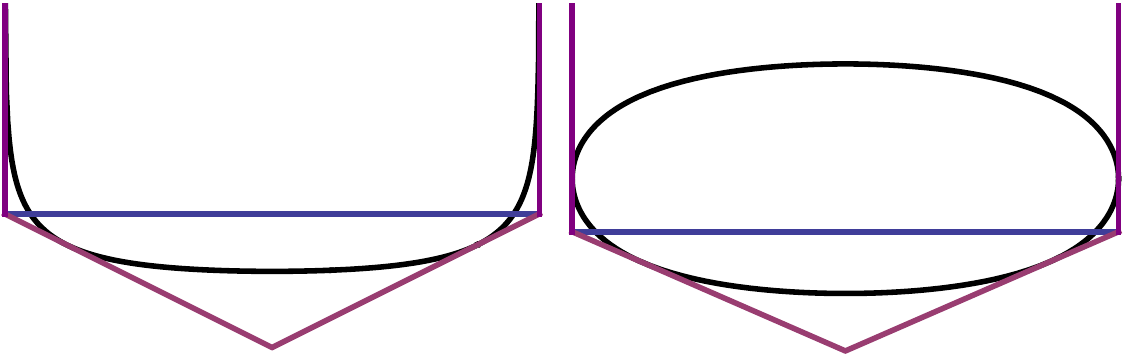}
\includegraphics[width=6.9cm]{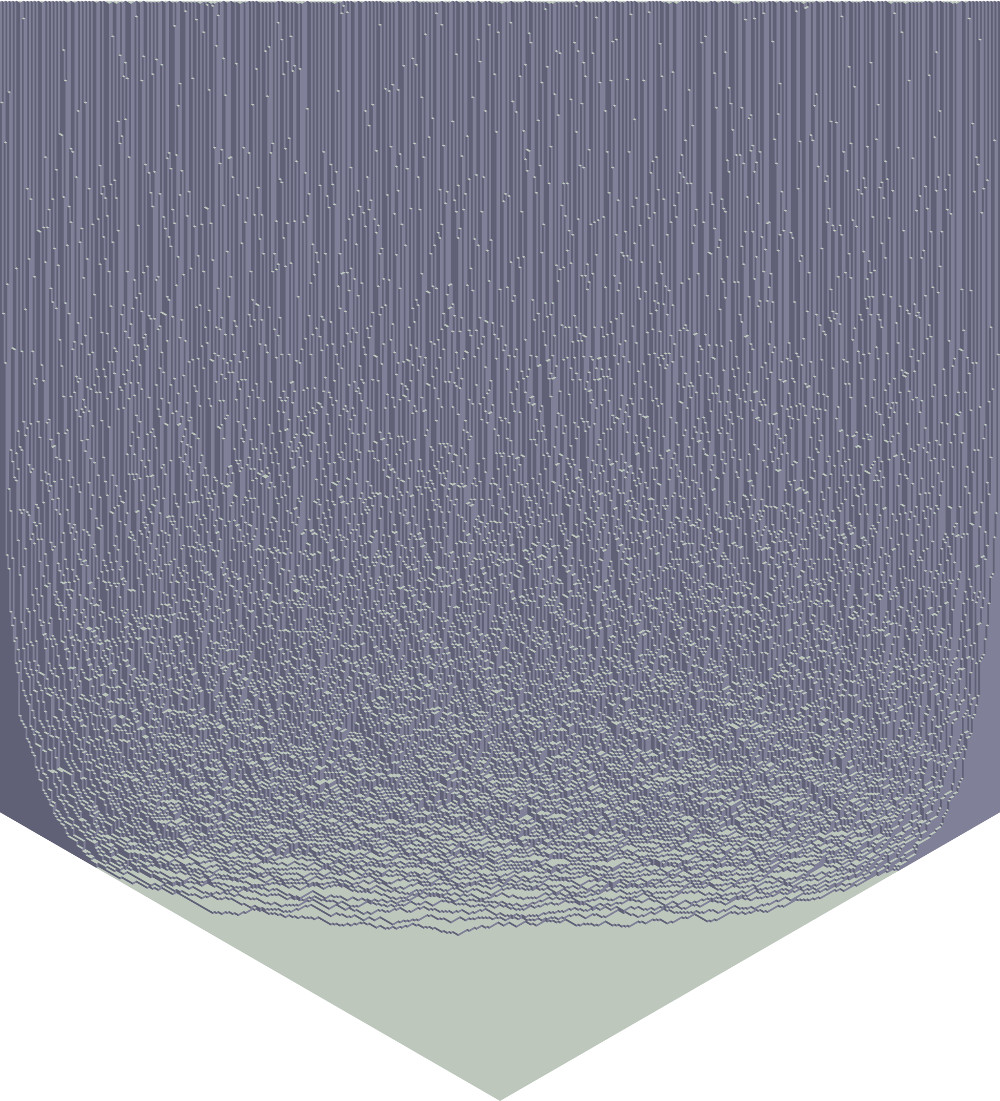}
\includegraphics[width=6.9cm]{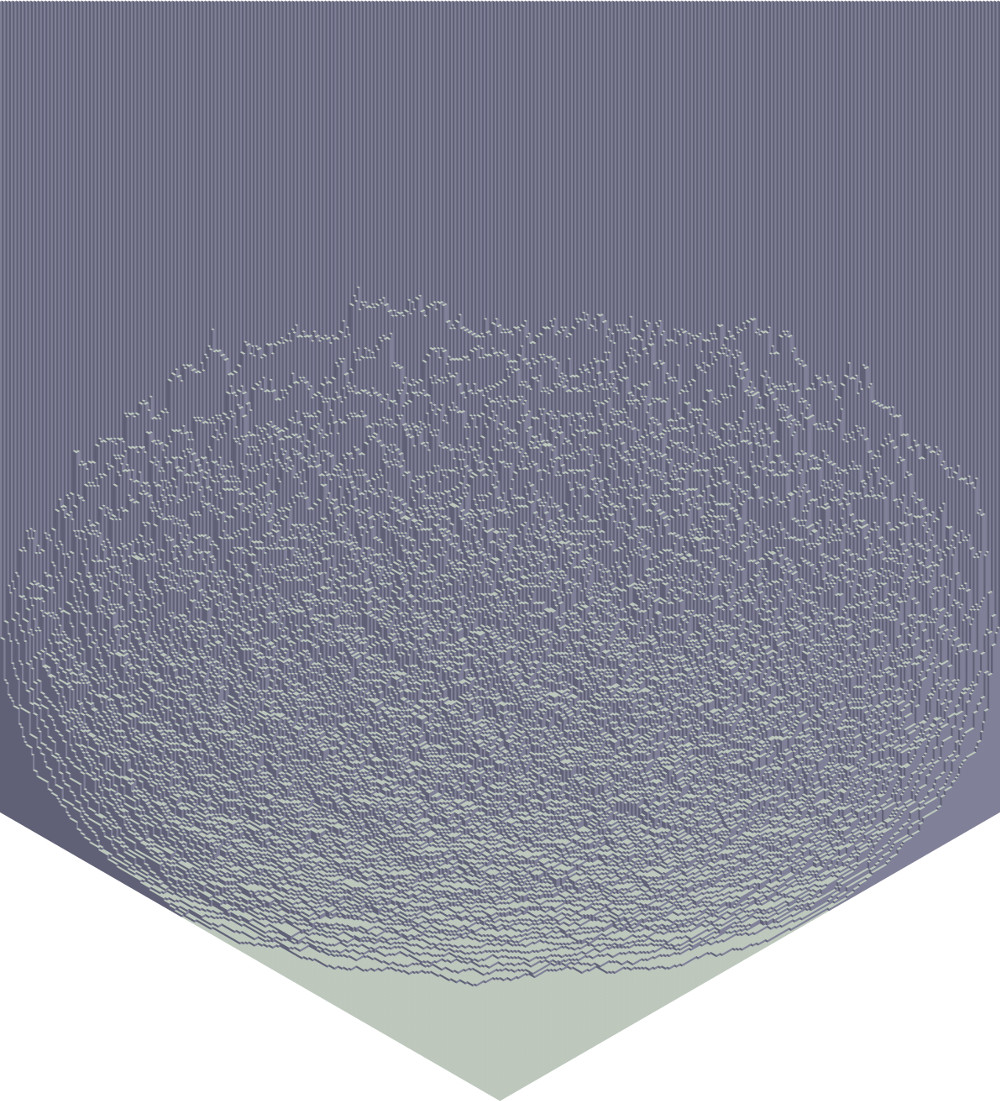}
\end{figure}

\end{itemize}

All exact samples were generated using Borodin's Schur Dynamics \cite{BorSchurDynam}.

\subsubsection{Turning points} 
\label{sec:turningPointsOverview}

In the bounded cases the system develops special points called turning points near the vertical boundaries, which we study in Section \ref{sec:turningPoints}. A turning point is a point where the disordered region meets two different frozen regions. In the case of homogeneous weights the system develops one turning point near each extreme. When the back wall is piecewise linear with lattice slopes these points were studied by Okounkov and Reshetikhin in \cite{OR3}, where it was proven that the local process at the turning point  is the same as the GUE minor process - the point process of the eigenvalues of the principal minors of a random $N\times N$ GUE matrix when $N\rightarrow\infty$. It was conjectured in \cite{OR3} that this holds universally, independently of boundary conditions. Results toward this were obtained in \cite{JohNordGUEminors} and \cite{GorinPanova}.

\begin{figure}[ht]
\caption{\label{fig:doubleTurningPoint-bottom} The bottom turning point on the right side of Figure \ref{fig:bddSample}. The picture is rotated 90 degrees clockwise.}
\includegraphics[width=14cm]{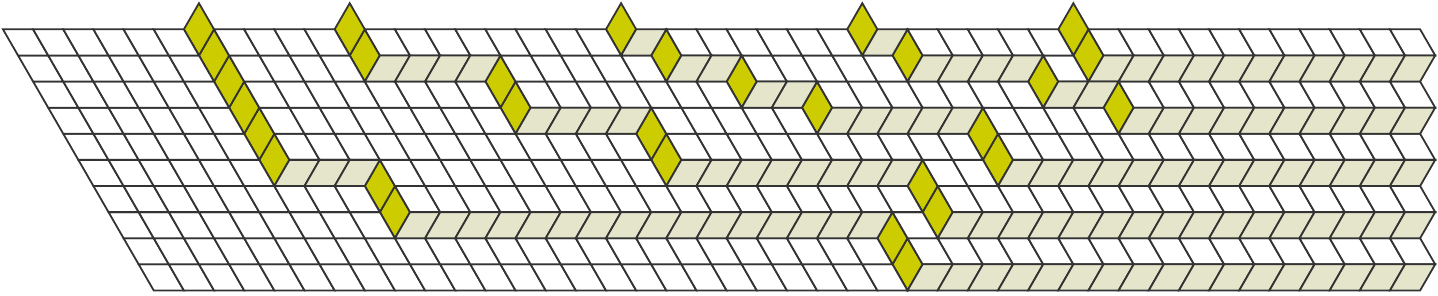}
\caption{\label{fig:doubleTurningPoint-top} The top turning point on the right side of Figure \ref{fig:bddSample}. The picture is rotated 90 degrees clockwise.}
\includegraphics[width=14cm]{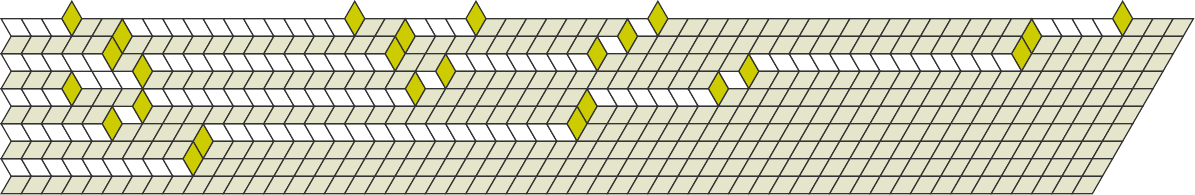}
\end{figure}

In the case of inhomogeneous weights the local process at turning points is not the GUE minor process. When $\alpha>1$, the turning point splits into two turning points separated by a frozen region where two types of lozenges coexist in a deterministic pattern - a so called semi-frozen region (see the right side of Figure \ref{fig:doubleTurningPoint-bottom} and left side of Figure \ref{fig:doubleTurningPoint-top}). Essentially at each turning point the system transitions from a wall with lattice slope to a wall with slope $0$. Whereas in the homogeneous case at a turning point at distance $n$ from the edge you see $n$ horizontal lozenges (corresponding to $n$ eigenvalues of an $n\times n$ matrix), in the inhomogeneous case you see $\lfloor \frac{n+1}2\rfloor$ or $\lfloor \frac n2\rfloor$ depending on whether it is the top or bottom turning point (see Figures \ref{fig:doubleTurningPoint-bottom} and \ref{fig:doubleTurningPoint-top}). Remarkably, if we only look at the points of distance of a fixed parity from the edge, we recover the GUE minor process. As a consequence at each turning point we see two non-trivially correlated copies of the GUE minor process.

It is noteworthy, that the GUE minor process we see when we look at the slices of only one parity implies interlacing of slices. This is not a result of a geometric constraint unlike the case of a single turning point found in the homogeneous case.

\subsubsection{An intermediate regime}
In Section \ref{sec:intermediateWeights} we consider the scaling limit of random skew plane partitions under the measure \eqref{eq:qinhomog} when the weights $q_t$ are given by
\begin{equation}
\label{eq:intermperiodicweights}
q_t=\left\{
\begin{array}{ll}
e^{-r+\gamma r^{1/2}},&t\text{ is even}\\
e^{-r-\gamma r^{1/2}},&t\text{ is odd}
\end{array}\right.,
\end{equation}
where $\gamma>0$ is an arbitrary constant. This is an intermediate regime between the homogeneous weights and the inhomogeneous weights given by \eqref{eq:periodicweights}. The macroscopic limit shape and correlations in the bulk are the same as in the homogeneous case, so periodicity disappears in the limit. However, the local point process at turning points is different from the homogeneous one. In particular, while we only have one turning point near each edge, we still do not have the GUE minor process, but rather a one-parameter deformation of it.

\subsection{Acknowledgements}
I am very grateful to Richard Kenyon for suggesting the study of periodic weights. I am very grateful to Richard Kenyon, Nicolai Reshetikhin, Paul Zinn-Justin, Andrei Okounkov and Vadim Gorin for many useful discussions on this subject. The project was started during the Random Spatial Processes program at MSRI and I am very grateful to MSRI with all its staff and the organizers of the special program for their excellent hospitality.

\section{Background and notation}
\subsection{Notation}

\begin{figure}[ht]
\caption{\label{fig:bLambda} The coordinates $u^\lambda_i$ and the graph of $b_\lambda(t)$ when $\lambda=\{3,2,2,1\}$ and $c=6,d=6$.}
\includegraphics[width=8cm]{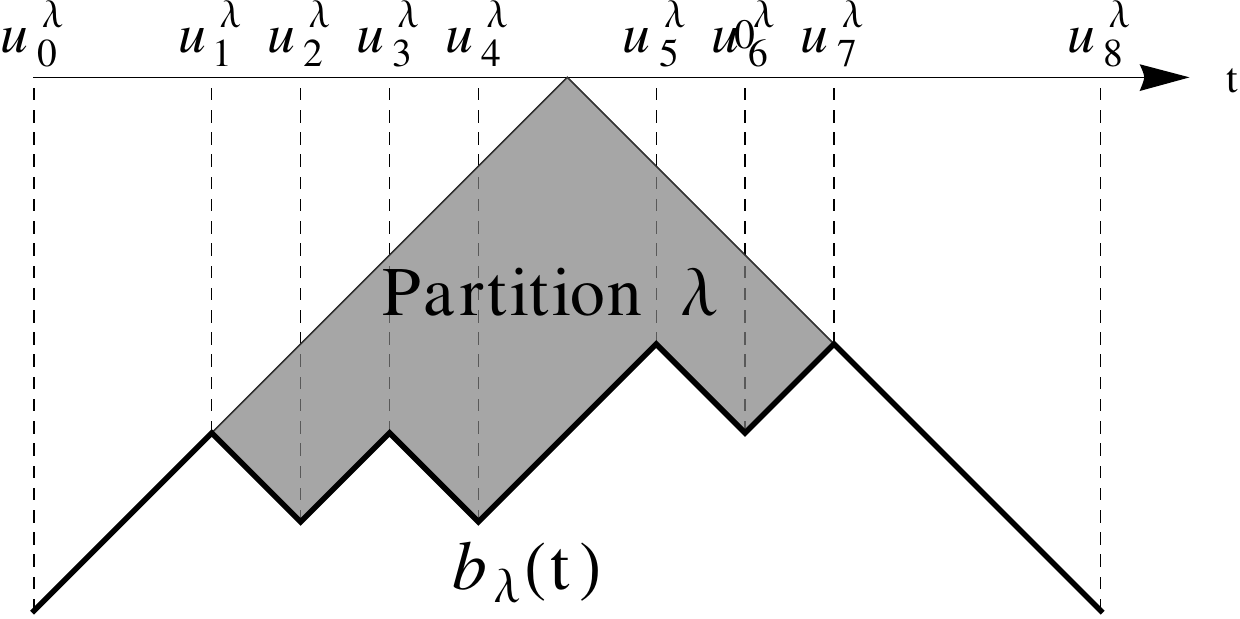}
\end{figure}

Let $u^\lambda_1<u^\lambda_2<\ldots<u^\lambda_{n-1}$ denote the horizontal coordinates of the corners on the boundary of the Young diagram $\lambda$. For convenience let $u^\lambda_0$ and $u^\lambda_n$ be the horizontal coordinates of respectively the left-most and right-most points of the back wall. We have $u^\lambda_0=-c$ and $u^\lambda_n=d$. Let
\begin{equation*}
 \tilde{u}^\lambda=\sum_{i=1}^{n-1}(-1)^i u^\lambda_i,
\end{equation*}
and for $t\in[c,d]$ define the piecewise linear function $b_\lambda(t)$ with slopes $\pm 1$ by
\begin{equation*}
b_\lambda(t)=\sum_{i=0}^{n} (-1)^{i}|t-\tilde{u}^\lambda-u^\lambda_i|+u^\lambda_0-u^\lambda_n.
\end{equation*}
The function $b_\lambda(t)$ gives the back wall of skew plane partitions with boundary $\lambda$ (see Figure \ref{fig:bLambda}). 

As can be seen from Figure \ref{fig:boxes}, skew plane partitions can be identified with tilings of regions of $\mathbb{R}^2$ with 3 types of rhombi, called lozenges. Scale the axes in such a way that the centers of horizontal tiles are on the lattice $\mathbb{Z}\times \frac 12 \mathbb{Z}$. Note, that the positions of the horizontal tiles completely determine the plane partition.

Given a subset $U=\{(t_1,h_1),\ldots,(t_n,h_n)\}\subset \mathbb{Z}\times \frac 12 \mathbb{Z}$, define the corresponding local correlation function $\rho_{\lambda,\bar{q}}(U)$ as the probability for a random tiling taken from the above probability space to have horizontal tiles centered at all positions $(t_i,h_i)_{i=1}^n$. 

Okounkov and Reshetikhin \cite{OR2} showed that for arbitrary $\lambda$ and $\bar{q}$, the point process of horizontal lozenges on $\mathbb{Z}\times\frac12\mathbb{Z}$ is determinantal with the following kernel:
\begin{theorem}[Theorem 2, part 3 \cite{OR2}] 
\label{thm:fin-corr2}
The correlation functions $\rho_{\lambda,\bar{q}}$ are determinants
\begin{equation*}
\rho_{\lambda,\bar{q}}(U)
=\det(K_{\lambda,\bar{q}}((t_i,h_i),(t_j,h_j)))_{1\leq i,j\leq n},
\end{equation*}
where the correlation kernel $K_{\lambda,\bar{q}}$ is given by the double integral
\begin{multline}
\label{eq:main-corr2}
K_{\lambda,\bar{q}}((t_1,h_1),(t_2,h_2))
=\\ \frac{1}{(2\pi \mathfrak{i})^2}
\int_{z\in C_z}\int_{w\in C_w}
\frac{\Phi_{b_\lambda,\bar{q}}(z,t_1)}{\Phi_{b_\lambda,\bar{q}}(w,t_2)}
\frac{1}{z-w}z^{-h_1+\frac 12 b_\lambda(t_1)+\frac 12}w^{h_2-\frac 12 b_\lambda(t_2)+\frac12}\frac{dz\ dw}{zw},
\end{multline}
where $b_\lambda(t)$ is the function giving the back wall corresponding to $\lambda$ as in Figure \ref{fig:bLambda}, 
\begin{align}
\nonumber\Phi_{b_\lambda,\bar{q}}(z,t)&=\frac{\Phi^-_{b_\lambda,\bar{q}}(z,t)}{\Phi^+_{b_\lambda,\bar{q}}(z,t)},\\
\label{eq:Phis}\Phi^+_{b_\lambda,\bar{q}}(z,t)&=\prod_{m>t, m\in D^+, m\in \mathbb{Z}+\frac 12}(1-zx_m^+),\\
\nonumber\Phi^-_{b_\lambda,\bar{q}}(z,t)&=\prod_{m<t, m\in D^-, m\in \mathbb{Z}+\frac 12}(1-z^{-1}x_m^-),
\end{align}
the parameters $x^\pm_m$ satisfy the conditions
\begin{align}
\label{eq:xs-qs}
\nonumber\frac{x^+_{m+1}}{x^+_m}&=q_{m+\frac 12},\ u^\lambda_{2i-1}<m<u^\lambda_{2i}-1,\ \text{or}\ m>u^\lambda_{n-1},
\nonumber\\x^+_{u^\lambda_{2i}-\frac 12}x^-_{u^\lambda_{2i}+\frac 12}&=q_{u^\lambda_{2i}}^{-1},
\nonumber\\x^-_{u^\lambda_{2i-1}-\frac 12}x^+_{u^\lambda_{2i-1}+\frac 12}&=q_{u^\lambda_{2i-1}},
\nonumber\\\frac{x^-_m}{x^-_{m+1}}&=q_{m+\frac 12},\ u^\lambda_{2i}<m<u^\lambda_{2i+1}-1,\ \text{or}\ m<u^\lambda_1,
\end{align}
$m\in D^{\pm}$ means the back wall at $t=m$, i.e. $b_\lambda(t)$ at $t=m$, has slope $\mp 1$, and $C_z$ (respectively $C_w$) is a simple positively oriented contour around 0 such that its interior contains none of the poles of $\Phi_{b_\lambda}(\cdot,t_1)$ (respectively all of the poles of $\Phi_{b_\lambda}(\cdot,t_2)^{-1}$). Moreover, if $t_1< t_2$, then $C_z$ is contained in the interior of $C_w$, and otherwise, $C_w$ is contained in the interior of $C_z$. 
\end{theorem}

Note, that the conditions \eqref{eq:xs-qs} can be rewritten as
\begin{align}
\label{eq:xs-prodqs}
x^-_m&=a^{-1}q_{u^\lambda_0+1}^{-1}\cdot\dots\cdot q_{m-\frac 12}^{-1},
\\\nonumber x^+_m&=aq_{u^\lambda_0+1}\cdot\dots\cdot q_{m-\frac 12},
\end{align}
where $a>0$ is an arbitrary parameter.

\subsection{The scaling limit}
\begin{figure}
\caption{\label{fig:backwall}The back wall}
\includegraphics[width=6cm]{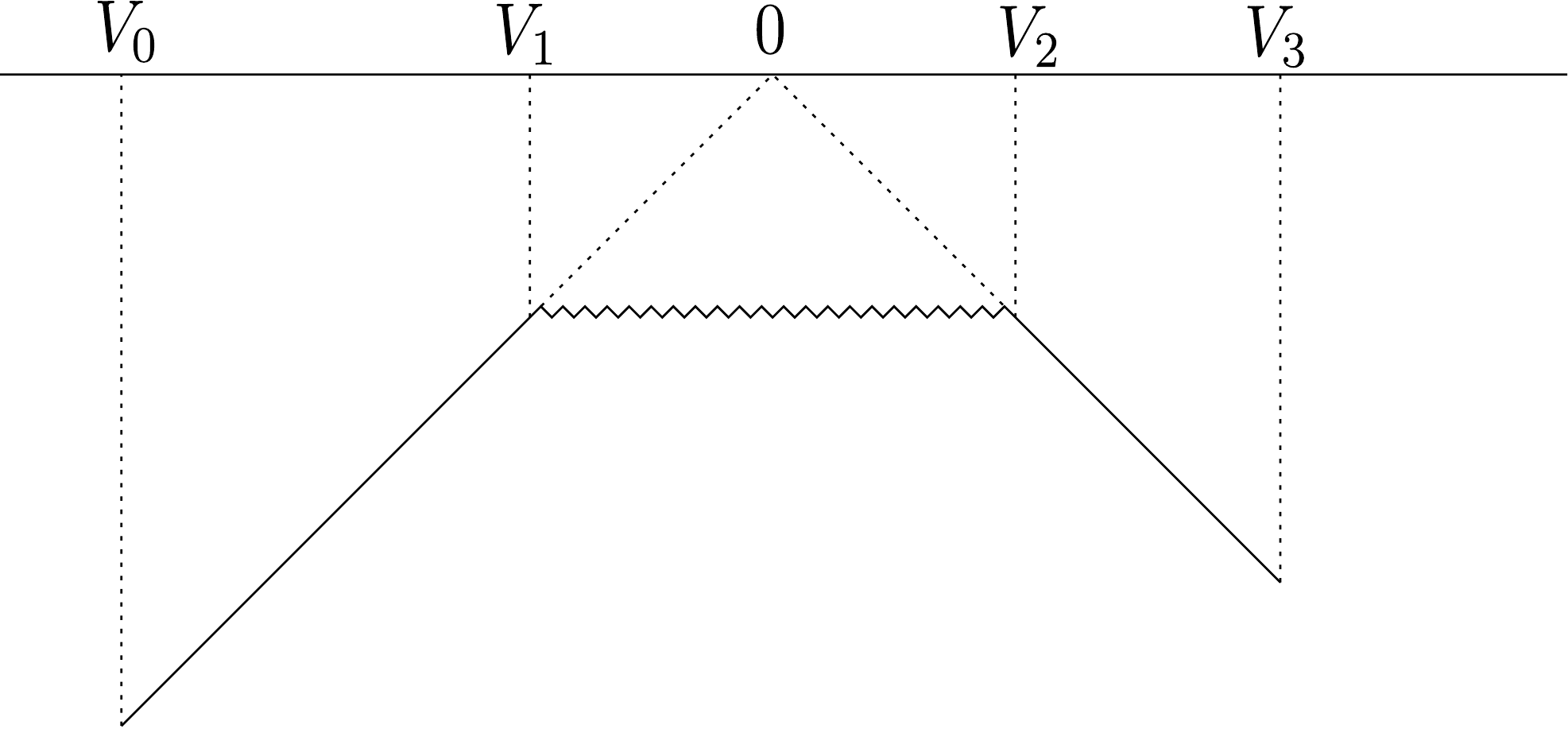}
\end{figure}

For $r>0$ let $\lambda_r$ be a staircase-shaped partition. Our goal is to study random skew plane partitions $\pi\in\Pi_{\lambda_r}^{c_r,d_r}$ under the measure \eqref{eq:qinhomog} with inhomogeneous weights given by \eqref{eq:periodicweights} in the limit when $r\rightarrow 0+$. Since the typical scale of such a random skew plane partition is $\frac 1r$, we scale plane partitions in all directions by $r$, and let $\tau=rt$, $\chi=rh$ be the rescaled coordinates. Let $B_{\lambda_r}(\tau)$ be the function giving the scaled back wall. We have $B_{\lambda_r}(\tau)=rb_{\lambda_r}(\tau/r)$. We assume that $\lim_{r\rightarrow 0} rc_r=-V_0$, $\lim_{r\rightarrow 0} rd_r=V_3$, and the scaled back walls given by $B_{\lambda_r}(\tau)$ converge point-wise and uniformly to the function $V(\tau)$ defined by
\begin{equation}
\label{eq:backwall}
V(\tau)=-\frac12|\tau-V_1|-\frac12|\tau-V_2|,
\end{equation}
where $V_1=-V_2$ (see Figure \ref{fig:backwall}).

\begin{figure}[ht]
\caption{\label{fig:stripWeight} Horizontal strip of width 1.}
\includegraphics[width=5cm]{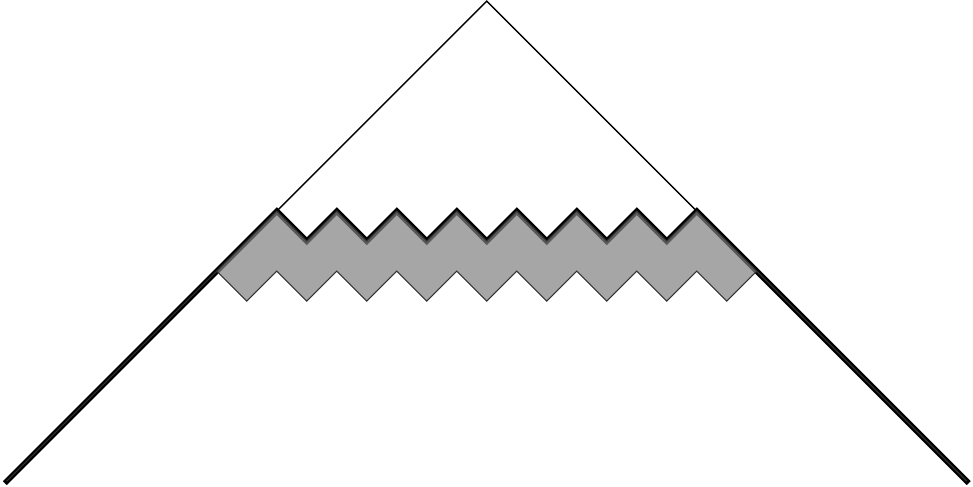}
\end{figure}

If $\alpha>1$, in order for the measure \eqref{eq:qinhomog} to be well defined, we must impose two constraints. First, the weight at inner corners must be less than $1$. Thus, we must have that $u^\lambda_1$ and $u^\lambda_{n-1}$ are even. For simplicity, we will also assume that $u^\lambda_0,u^\lambda_n$ are odd. Second, the length of the section with average slope $0$ should be large enough so that the total weight of the strip of width $1$ in this section is less than $1$ (see Figure \ref{fig:stripWeight}). Since we are scaling by $r$, there are $(V_2-V_1)/r$ microscopic linear sections between $V_1$ and $V_2$. Hence, the weight of the strip in question is 
\begin{equation*}
\alpha(e^{-r})^{\frac{V_2-V_1}r}.
\end{equation*}
Thus, we must have
\begin{equation}
\label{eq:alphaRestriction}
e^{-(V_2-V_1)}\alpha<1.
\end{equation}

\section{The correlation kernel in the bulk}
\label{sec:bulkCorrKer}

\subsection{Computation of the exponentially leading term of the correlation kernel}

In order to understand the asymptotic local correlations near a given macroscopic point $(\tau,\chi)$, we need to study the limit of the correlation kernel \eqref{eq:main-corr2} when $\lim_{r\rightarrow 0}rt_1=\lim_{r\rightarrow 0}rt_2=\tau$, $\lim_{r\rightarrow 0}rh_1=\lim_{r\rightarrow 0}rh_2=\chi$, $\Delta t=t_1-t_2$ and $\Delta h=h_1-h_2$ are constants, and $t_1, t_2$ have fixed parity, independent of $r$. To do this, we apply the saddle point method to the integral representation of the correlation kernel \eqref{eq:main-corr2}. The first step is to compute the exponentially leading term of the integrand.

From \eqref{eq:Phis} we have
\begin{align*}
-r\ln\Phi_+(z,t_i)
=&-\sum_{m>t_i}r\frac 12(1-b'_{\lambda_r}(m))\ln(1-e^{-r(m-\frac 12-u^\lambda_0)}\alpha^{-p_m}a_r z)
\\=&-\sum_{m>t_i}r\frac 12(1-B'_{\lambda_r}(m))\ln(1-e^{-rm}\alpha^{-p_m}\tilde{a}_r z),
\end{align*}
where
\begin{equation*}
p_m=\left\{
\begin{array}{ll}
0,&m-\frac 12 - u^\lambda_0 \text{ is even}
\\1,&m-\frac 12 - u^\lambda_0 \text{ is odd}
\end{array}
\right.,
\end{equation*}
and $\tilde{a}_r=a_re^{r(\frac12+u^\lambda_0)}$.
Making a change of variable, we obtain
\begin{align*}
-r\ln\Phi_+(z,t_i)
=&-\frac 12\sum_{M\in(\tau,\infty)\cap r(2\mathbb{Z}+\frac 12)}r(1-B'_{\lambda_r}(M))\ln(1-e^{- M}\alpha^{-1}\tilde{a}_rz)
\\&-\frac 12\sum_{M\in(\tau,\infty)\cap r(2\mathbb{Z}+\frac32)}r(1-B'_{\lambda_r}(M))\ln(1-e^{- M}\tilde{a}_rz),
\end{align*}
where $\tau=rt_i$.

Since $B_{\lambda_r}(M)$ converges to $V(M)$, we have
\begin{align*}
\lim_{k\rightarrow\infty}-r\ln\Phi_+(z,t_i)
=&-\frac12\int\limits_{(\tau,\infty)\cap(V_2,V_3)}\ln(1-e^{- M}\tilde{a}\alpha^{-1}z)+\ln(1-e^{- M}\tilde{a}z)dM
\\&-\frac12\int\limits_{(\tau,\infty)\cap(V_1,V_2)}\ln(1-e^{- M}\tilde{a}\alpha^{-1}z)dM.
\end{align*}

Similar computations for $\Phi_-(z,t_i)$ and setting the arbitrary parameter $\tilde{a}_r$ to $\alpha^{\frac 12}$, 
give us

\begin{equation}
\label{eq:leadAsympCorrKer}
K_{\lambda_r,\bar{q}_r}((t_1,h_1),(t_2,h_2))
=\frac{1}{(2\pi \mathfrak{i})^2}
\int_{z\in C_z}\int_{w\in C_w}
e^{\frac{S_{\tau,\chi}(z)-S_{\tau,\chi}(w)}r+O(1)}\frac{dz\ dw}{z-w},
\end{equation}
where
%
%
\begin{align}
\label{eq:S}
S_{\tau,\chi}(z)
=&\frac12\int\limits_{V_0}^{\min(V_1,\tau)}\ln(1-e^{ M}\alpha^{\frac 12} z^{-1})+\ln(1-e^{ M}\alpha^{-\frac 12}z^{-1})dM
\\\nonumber&+\frac12\int\limits_{\min(V_1,\tau)}^{\min(V_2,\tau)}\ln(1-e^{ M}\alpha^{-\frac 12}z^{-1})dM
\\\nonumber&-\frac12\int\limits_{\max(\tau,V_2)}^{V_3}\ln(1-e^{- M}\alpha^{-\frac 12}z)+\ln(1-e^{- M}\alpha^{\frac 12}z)dM
\\\nonumber&-\frac12\int\limits_{\max(\tau,V_1)}^{\max(\tau,V_2)}\ln(1-e^{- M}\alpha^{-\frac 12}z)dM
-(\chi-\frac 12 V(\tau))\ln z.
\end{align}


\subsection{Critical points}
To apply the saddle point method we need to study the critical points of the function $S_{\tau,\chi}$. To simplify formulas, we set $v=V_3=-V_0$, and $u=V_2=-V_1$. Differentiating the formula for $S_{\tau,\chi}(z)$ and using \eqref{eq:backwall} we obtain
\begin{align}
\label{eq:spBounded}
 z\frac{d}{dz}S_{\tau,\chi}(z)
=&-\frac12\ln\left(\frac{z-e^{\min(-u,\tau)}\alpha^{\frac12}}{z-e^{- v}\alpha^{\frac12}}\right)
-\frac12\ln\left(\frac{z-e^{\min(u,\tau)}\alpha^{-\frac12}}{z-e^{- v}\alpha^{-\frac12}}\right)
\\\nonumber&+\frac12\ln\left(\frac{z-e^{ v}\alpha^{\frac12}}{z-e^{\max(-u,\tau)}\alpha^{\frac12}}\right)
+\frac12\ln\left(\frac{z-e^{ v}\alpha^{-\frac12}}{z-e^{\max(u,\tau)}\alpha^{-\frac12}}\right)
\\\nonumber&-(\chi-\frac12 \tau+v).
\end{align}

Let $z_{\tau,\chi}$ denote a critical point of $S_{\tau,\chi}(z)$. Define
\begin{align*}
\mathfrak{U}:
=&(e^{- v}\alpha^{\frac12},e^{\min(-u,\tau)}\alpha^{\frac12})
\cup(e^{- v}\alpha^{-\frac12},e^{\min(u,\tau)}\alpha^{-\frac12})
\\&\cup(e^{\max(-u,\tau)}\alpha^{\frac12},e^{ v}\alpha^{\frac12})
\cup(e^{\max(u,\tau)}\alpha^{-\frac12},e^{ v}\alpha^{-\frac12}).
\end{align*}
It is easy to check that the critical points are not near the set $\mathfrak{U}$, in the sense that for any $(\tau,\chi)$ we have 
\begin{equation*}
\min_{x\in\mathfrak{U}}|z_{\tau,\chi}-x|>0. 
\end{equation*}


\subsubsection{The number of complex critical points everywhere}
We show that for any pair $(\tau,\chi)$, $S_{\tau,\chi}$ has at most one non-real complex conjugate pair of critical points. 

Let
\begin{equation*}
P_1(z)=e^{\chi-\frac12\tau+v}
(z-e^{- v}\alpha^{-\frac12})
(z-e^{- v}\alpha^{\frac12})
(z-e^{ v}\alpha^{-\frac12})
(z-e^{ v}\alpha^{\frac12})
\end{equation*}
and
\begin{equation*}
P_2(z)=
(z-e^{ \min(-u,\tau)}\alpha^{\frac12})
(z-e^{ \min(u,\tau)}\alpha^{-\frac12})
(z-e^{ \max(-u,\tau)}\alpha^{\frac12})
(z-e^{ \max(u,\tau)}\alpha^{-\frac12}).
\end{equation*}

Exponentiating $2 z\frac d{dz}S_{\tau,\chi}(z)$ we see that $P_1(z_{\tau,\chi})-P_2(z_{\tau,\chi})=0$. It follows from $-v<u,\tau<v$ and \eqref{eq:alphaRestriction} that
\begin{equation*}
e^{- v}\alpha^{-\frac12}
<e^{- v}\alpha^{\frac12}
<e^{ \max(u,\tau)}\alpha^{-\frac12}
<e^{ v}\alpha^{-\frac12}
<e^{ v}\alpha^{\frac12}
\end{equation*}
and
\begin{equation*}
e^{- v}\alpha^{-\frac12}
<e^{ \min(-u,\tau)}\alpha^{\frac12},
e^{ \min(u,\tau)}\alpha^{-\frac12},
e^{ \max(-u,\tau)}\alpha^{\frac12},
e^{ \max(u,\tau)}\alpha^{-\frac12}
<e^{ v}\alpha^{\frac12}.
\end{equation*}
Thus,
\begin{align*}
P_1(e^{- v}\alpha^{-\frac12})-P_2(e^{- v}\alpha^{-\frac12})<0,\\
P_1(e^{ v}\alpha^{\frac12})-P_2(e^{ v}\alpha^{\frac12})>0,
\end{align*}
and
\begin{equation*}
P_1(e^{ \max(u,\tau)}\alpha^{-\frac12})-P_2(e^{ \max(u,\tau)}\alpha^{-\frac12})<0.
\end{equation*}
By the intermediate value theorem $P_1(z)-P_2(z)$ has at least two distinct real roots. Since $P_1(z)-P_2(z)$ is a degree $4$ polynomial in $z$, the number of non-real complex roots is either $0$ or $2$. Equivalently, the number of non-real complex critical points of $S_{\tau,\chi}$ is either $0$ or $2$.

\subsubsection{Critical points when $|\chi|\gg 1$}
Fix $\tau$. If $|\chi|\rightarrow \infty$, then looking at the real part of \eqref{eq:spBounded} it is easy to see that $|z_{\tau,\chi}|\rightarrow e^{ v}\alpha^{\pm \frac12}$ or $e^{- v}\alpha^{\pm \frac12}$ if $\chi<0$, and $|z_{\tau,\chi}|\rightarrow e^{\min(-u,\tau)}\alpha^{\frac 12}$, $e^{\max(-u,\tau)}\alpha^{\frac 12}$, $e^{\min(u,\tau)}\alpha^{-\frac 12}$ or $e^{\max(u,\tau)}\alpha^{-\frac 12}$ if $\chi>0$. In all these cases a direct application of \cite[Lemma 2.3]{M} shows that if $\Im\left( z\frac{d}{dz}S_{\tau,\chi}(z_{\tau,\chi})\right)= 0$, then $z_{\tau,\chi}$ must be real unless $\tau=\pm u$, $\chi>0$ and $u\neq v$. In the latter case $S_{\tau,\chi}$ has exactly one pair of complex conjugate critical points, the asymptotically leading term of which can be explicitly computed to give
\begin{multline*}
z_{\tau,\chi}=e^{\tau}\alpha^{\mp\frac12}+\mathfrak{i} 
\frac{
(e^{\tau}\alpha^{\mp\frac12}-e^{-v}\alpha^{\frac12})^{\frac12}
(e^{\tau}\alpha^{\mp\frac12}-e^{-v}\alpha^{-\frac12})^{\frac12}}
{(e^{\tau}\alpha^{\mp\frac12}-e^{-u}\alpha^{\pm\frac12})^{\frac12}
}
\\\times\frac{(e^{v}\alpha^{\frac12}-e^{\tau}\alpha^{\mp\frac12})^{\frac12}
(e^{v}\alpha^{-\frac12}-e^{\tau}\alpha^{\mp\frac12})^{\frac12}}
{(e^{\tau}\alpha^{\mp\frac12}-e^{u}\alpha^{\pm\frac12})^{\frac12}}e^{-\chi+\frac12\tau-v+O(e^{-\chi})}
\end{multline*}
or its conjugate, where the top sign is chosen when $\tau=u$ and the bottom sign when $\tau=-u$.

\subsection{The correlation kernel in the bulk}
Suppose $(\tau,\chi)$ is such that $S_{\tau,\chi}$ has a pair of complex conjugate critical points. Let the critical points be $z_{\tau,\chi}$ and $\bar{z}_{\tau,\chi}$, with $\Im z_{\tau,\chi}>0$. Following the saddle point method, we deform the contours of integration in \eqref{eq:main-corr2} to new contours $C'_z$, $C'_w$ as follows:
\begin{itemize}
\item the contours $C'_z$, $C'_w$ pass through the critical points $z_{\tau,\chi}$, $\bar{z}_{\tau,\chi}$,
\item along the contour $C'_z$ we have $\Re S_{\tau,\chi}(z)\leq \Re S_{\tau,\chi}(z_{\tau,\chi})$, with equality only at the critical points,
\item along the contour $C'_w$ we have $\Re S_{\tau,\chi}(w)\geq \Re S_{\tau,\chi}(z_{\tau,\chi})$, with equality only at the critical points.
\end{itemize}
During this deformation the contours cross each other along an arc connecting the conjugate critical points, and we pick up residues along the arc. The arc will cross the real axis at a positive point if $t_1\geq t_2$ and at a negative point otherwise. Since along the new contours we have 
\begin{equation*}
\Re S_{\tau,\chi}(z)\leq \Re S_{\tau,\chi}(z_{\tau,\chi})\leq \Re S_{\tau,\chi}(w),
\end{equation*}
the contours cross transversally at the critical points and the leading term of the integrand is $e^{\frac{S_{\tau,\chi}(z)-S_{\tau,\chi}(w)}r}$, 
as $r\rightarrow 0$, the contribution of the double integral along the new contours is exponentially small and the main contribution comes from the residue term. Thus, we have
\begin{multline*}
\lim_{r\rightarrow 0}K_{\lambda_r,\bar{q}_r}((t_1,h_1),(t_2,h_2))
=\\\lim_{r\rightarrow 0} \frac{1}{2\pi \mathfrak{i}}
\int_{\bar{z}_{\tau,\chi}}^{z_{\tau,\chi}}
\frac{\Phi_{b_{\lambda_r},\bar{q}_r}(z,t_1)}{\Phi_{b_{\lambda_r},\bar{q}_r}(z,t_2)}
z^{-\Delta h+\frac 12 (b_{\lambda_r}(t_1)-b_{\lambda_r}(t_2))-1} dz.
\end{multline*}

For $e\in\{0,1\}$ denote 
\begin{multline*}
m_e(t_1,t_2)=\sign(\Delta t)
\\\times
\left|\left\{m\in D^-:\min(t_1,t_2)<m<\max(t_1,t_2)
\&m-\frac12\in2\mathbb{Z}+e\right\}\right|,
\end{multline*}
and let $m(t_1,t_2)=m_0(t_1,t_2)+m_1(t_1,t_2)$.
Using \eqref{eq:Phis}, we obtain
\begin{multline}
\label{eq:corr-bulk}
\lim_{r\rightarrow 0}K_{\lambda_r,\bar{q}_r}((t_1,h_1),(t_2,h_2))
=\lim_{r\rightarrow 0} 
(-1)^{m(t_1,t_2)}(\alpha^{-1/2}e^\tau)^{m_0(t_1,t_2)}(\alpha^{1/2}e^\tau)^{m_1(t_1,t_2)}
\\\times\frac{1}{2\pi \mathfrak{i}}
\int_{\bar{z}_{\tau,\chi}}^{z_{\tau,\chi}}
z^{-\Delta h-\frac 12 \Delta t-1} \prod_{\stackrel{\min(t_1,t_2)<m<\max(t_1,t_2)}{m-\frac12\in\mathbb{Z}}}(1-zx_m^+)^{\sign \Delta t}
dz.
\end{multline}

Since the correlations are given by determinants and $m_e(t_1,t_2)$ is of the form $f(t_1)-f(t_2)$, the terms outside of the integral will cancel when taking determinants. Removing those terms we obtain the following theorem.

\begin{theorem}
The correlation functions of the system near a point $(\tau, \chi)$ in the bulk are given by
\begin{equation*}
K^{\alpha}_{\chi,\tau}(t_1,t_2,\Delta h) = \int_\gamma (1- e^{-\tau}\alpha^{\frac12} z)^{\frac{\Delta t+c}2} (1- e^{-\tau}\alpha^{-\frac12} z)^{\frac{\Delta t-c}2} z^{-\Delta h -\frac{\Delta t}{2}} \frac{dz}{2\pi\mathfrak{i} z},
\end{equation*}
where
\begin{equation*}
c=\left\{
\begin{array}{rl}
1,&\Delta t\text{ is odd and }t_1\text{is even}\\
0,&\Delta t\text{ is even}\\
-1,&\text{otherwise}
\end{array}
\right.,
\end{equation*}
the integration contour connects the two non-real critical points of $S_{\tau,\chi}(z)$, passing through the real line in the interval $(0,e^{\tau}\alpha^{-\frac 12})$ if $\Delta(t)\geq 0$ and through $(-\infty, 0)$ otherwise.
\end{theorem} 
\begin{remark}
The local point process in the bulk is not $\mathbb{Z}\times\mathbb{Z}$ invariant as in the homogeneous case, but rather $2\mathbb{Z}\times \mathbb{Z}$ translation invariant.
\end{remark}

\section{The frozen boundary}
\label{sec:frozenBoundary}

The region of the $(\tau,\chi)$ plane consisting of the points $(\tau,\chi)$ where $S_{\tau,\chi}$ has complex conjugate critical points is called the disordered region. Suppose $(\tau,\chi)$ is in the complement of the closure of the disordered region. As we deform the contours of integration in \eqref{eq:main-corr2} following the saddle point method, the contours either do not cross at all, in which case the correlation kernel converges to $0$, or cross along a closed curve winding once around the origin, in which case the correlation kernel converges to $1$. Thus, at such points $(\tau,\chi)$ horizontal lozenges appear with probability $0$ or $1$ and we have a frozen region. 

The boundary of the disordered region is called the frozen boundary. It consists of points $(\tau,\chi)$ such that $S_{\tau,\chi}$ has double real critical points, i.e. points $(\tau,\chi)$ such that there exist $z\in\mathbb{R}$ satisfying $S'_{\tau,\chi}(z)=S''_{\tau,\chi}(z)=0$. We study this curve in three different regimes.

%

\subsection{Infinite floor}
Consider skew plane partitions with an unbounded floor. In our notation this means $V_0=-\infty$, $V_3=\infty$ and $u=V_2=-V_1$.

We are interested in the points $(\tau,\chi)$ where $S_{\tau,\chi}$ has double real critical points. We will show that for every $z\in\mathbb{R}\backslash\{0\}$, there is a unique pair $(\tau,\chi)$ for which it is a double real critical point. Thus, the set of points $(\tau,\chi)$ with double real critical points can be parametrized by $z$. We will write $(\tau(z),\chi(z))$ for this curve.

Differentiating \eqref{eq:spBounded} with $v=\infty$ we obtain
\begin{align}
\label{eq:sppInfinite}
\frac{d}{dz}\left( z\frac{d}{dz} S_{\tau,\chi}(z)\right)
=&\frac1z-\frac12\frac1{z-e^{\tau}\alpha^{\frac12}}
-\frac12\frac1{z-e^{\tau}\alpha^{-\frac12}}
\\\nonumber&-\frac12\frac1{z-e^{- u}\alpha^{\frac12}}
-\frac12\frac1{z-e^{ u}\alpha^{-\frac12}}.
\end{align}

We can solve the equation
\begin{equation}
\label{eq:sppInfiniteIs0}
\frac{d}{dz}\left( z\frac{d}{dz} S_{\tau,\chi}(z)\right)=0
\end{equation}
for $\tau$ in terms of $z$. Once we have $z$ and $\tau$, $\chi$ can be determined uniquely from $ z\frac{d}{dz} S_{\tau,\chi}(z)=0$.

Setting
$$\mathfrak{A}=\alpha^{\frac 12}+\alpha^{-\frac 12}$$
and
$$\mathfrak{B}=e^{- u}\alpha^{\frac12}+e^{ u}\alpha^{-\frac12},$$ equation \eqref{eq:sppInfiniteIs0} is equivalent to
\begin{equation}
\label{eq:sppInfiniteQuadratic}
e^{2\tau}(\mathfrak{B}z-2)-e^{\tau}\mathfrak{A}z(z^2-1)+z^3(2z-\mathfrak{B})=0,
\end{equation}
which is quadratic in $e^{\tau}$. Let $\tau^\pm(z)$ be the solutions
\begin{equation}
\label{eq:tauPMz}
e^{\tau^\pm(z)}=\frac{\mathfrak{A}z(z^2-1)\pm\sqrt{\mathfrak{A}^2(z^3-z)^2-4z^3(\mathfrak{B}z-2)(2z-\mathfrak{B})}}{2(\mathfrak{B}z-2)}.
\end{equation}
We show that for any $z\in\mathbb{R}$ only one of the solutions leads to a real pair $(\tau,\chi)$.

If $z<0$, then the leading coefficient of \eqref{eq:sppInfiniteQuadratic} is negative while the free coefficient is positive, whence the roots have opposite signs. Since $\tau\in\mathbb{R}$, we must have $\tau(z)=\tau^-(z)$.

Suppose $z>0$. First, notice that if $\frac2{\mathfrak{B}}<z<\frac{\mathfrak{B}}2$, then the leading coefficient of \eqref{eq:sppInfiniteQuadratic} is positive while the free coefficient is negative. Thus, $e^{\tau^-(z)}<0<e^{\tau^+(z)}$ and $\tau^-(z)\notin\mathbb{R}$. 

Recall, that if $z$ is a critical point, we must have $z\notin \mathfrak{U}$. We will show, that this condition is satisfied by exactly one of the solutions $\tau^\pm(z)$, namely $\tau^+(z)$. To show this, it is enough to show the following:
\begin{align}
\label{eq:forbiddenInfinite1}
0<z<e^{- u}\alpha^{\frac12}
&\Rightarrow e^{\tau^+(z)}\alpha^{\frac12}<
z<e^{\tau^-(z)}\alpha^{\frac12},
\\\label{eq:forbiddenInfinite2}
e^{- u}\alpha^{\frac12}<z<\frac{2}{\mathfrak{B}}
&\Rightarrow e^{\tau^+(z)}\alpha^{-\frac12}<
z<e^{\tau^-(z)}\alpha^{-\frac12},
\\\label{eq:forbiddenInfinite3}
\frac{\mathfrak{B}}2<z<e^{ u}\alpha^{-\frac12}
&\Rightarrow e^{\tau^-(z)}\alpha^{\frac12}<
z<e^{\tau^+(z)}\alpha^{\frac12},
\\\label{eq:forbiddenInfinite4}
e^{ u}\alpha^{-\frac12}<z
&\Rightarrow e^{\tau^-(z)}\alpha^{-\frac12}<
z<e^{\tau^+(z)}\alpha^{-\frac12}.
\end{align}
We will show \eqref{eq:forbiddenInfinite1}. The remaining three can be established similarly. If $0<z<e^{- u}\alpha^{\frac12}$, then the leading coefficient in \eqref{eq:sppInfiniteQuadratic} is negative, so showing \eqref{eq:forbiddenInfinite1} is equivalent to showing that the left-hand side of \eqref{eq:sppInfiniteQuadratic} evaluated at $e^{\tau}=z\alpha^{-\frac12}$ is positive. Substituting $e^{\tau}=z\alpha^{-\frac12}$ in \eqref{eq:sppInfiniteQuadratic} we obtain
\begin{equation*}
(1-\alpha^{-1})z^2(z-e^{- u}\alpha^{\frac12})(z-e^{ u}\alpha^{-\frac12}),
\end{equation*}
which is negative, since $z<e^{- u}\alpha^{\frac 12}<e^{ u}\alpha^{-\frac 12}$ and $\alpha>1$.

Thus, we obtain that for any $z\in\mathbb{R}$, there is a unique pair $(\tau(z),\chi(z))$ such that $z$ is a double real critical point for $S_{\tau,\chi}$. This curve is the frozen boundary. From \eqref{eq:tauPMz} it is easy to read general features of the frozen boundary such as the appearance of the tentacles, as
\begin{equation*}
\lim_{z\rightarrow 0}\tau(z)=-\infty
,\qquad \lim_{z\rightarrow \pm\infty}\tau(z)=\infty
,\qquad \lim_{z\rightarrow (e^{- u}\alpha^{\frac12})^\pm}\tau(z)=-u^\pm,
\end{equation*}
and
\begin{equation*}
\lim_{z\rightarrow (e^{ u}\alpha^{-\frac12})^\pm}\tau(z)=u^\pm
,\qquad \lim_{z\rightarrow (e^{ u}\alpha^{-\frac12})^\pm}\chi(z)=\infty
,\qquad \lim_{z\rightarrow (e^{- u}\alpha^{\frac12})^\pm}\chi(z)=\infty,
\end{equation*}
where for example by an expression like $\lim_{z\rightarrow a^+}\tau=b^-$ we mean $\tau$ approaches $b$ from below when $z$ approaches $a$ from above.


\subsection{Bounded floor}
Let us examine the frozen boundary near the edges $\tau=\pm v$ when the floor is bounded. Suppose $(\tau,\chi)$ is a point of the frozen boundary and $z_{\tau,\chi}$ is the double real critical point of $S_{\tau,\chi}$. Differentiating \eqref{eq:spBounded} we obtain
\begin{equation*}
Q_1(z_{\tau,\chi})+Q_2(z_{\tau,\chi})=0,
\end{equation*}
where
\begin{equation*}
Q_1(z)=\frac{1}{z-e^v\alpha^{\frac12}}
-\frac{1}{z-e^{\tau}\alpha^{\frac12}}
+\frac{1}{z-e^v\alpha^{-\frac12}}
-\frac{1}{z-e^{\tau}\alpha^{-\frac12}}
\end{equation*}
and
\begin{equation*}
Q_2(z)=
\frac{1}{z-e^{-v}\alpha^{-\frac12}}
+\frac{1}{z-e^{-v}\alpha^{\frac12}}
-\frac{1}{z-e^{-u}\alpha^{\frac12}}
-\frac{1}{z-e^u\alpha^{-\frac12}}.
\end{equation*}
Since 
\begin{equation*}
e^{-v}\alpha^{-\frac12}
<e^{-v}\alpha^{\frac12}
<e^{-u}\alpha^{\frac12}
<e^u\alpha^{-\frac12},
\end{equation*}
the equation $Q_2(z)=0$ has two real solutions, both in the interval $(e^{-v}\alpha^{-\frac12},e^u\alpha^{-\frac12})$. However, if $\tau>u$, then $(e^{-v}\alpha^{-\frac12},e^u\alpha^{-\frac12})\subset \mathfrak{U}$ and it follows from $z_{\tau,\chi}\notin\mathfrak{U}$ that $\lim_{\tau\rightarrow v}Q_1(z_{\tau,\chi})\neq 0$. Thus
\begin{equation*}
\lim_{\tau\rightarrow v}z_{\tau,\chi}=e^v\alpha^{\frac12}
\text{ or }
\lim_{\tau\rightarrow v}z_{\tau,\chi}=e^v\alpha^{-\frac12}.
\end{equation*}
Similarly,
\begin{equation*}
\lim_{\tau\rightarrow -v}z_{\tau,\chi}=e^{-v}\alpha^{\frac12}
\text{ or }
\lim_{\tau\rightarrow -v}z_{\tau,\chi}=e^{-v}\alpha^{-\frac12}.
\end{equation*}
It follows that there are two turning points near each of the extremes $\tau=\pm v$. The vertical coordinates of the turning points are
\begin{equation}
\label{eq:chiBottom}
\chi_{bottom}=-\frac{v}2-\frac12\ln
\frac{(e^v-e^{-u})(\alpha e^v-e^u)}{(e^v-e^{-v})(\alpha e^v-e^{-v})}
\end{equation}
corresponding to $z_{\tau,\chi}=e^\tau\alpha^{\frac12}$ and
\begin{equation}
\label{eq:chiTop}
\chi_{top}=-\frac{v}2-\frac12\ln
\frac{(e^v-e^{-u}\alpha)(e^v-e^u)}{(e^v-e^{-v}\alpha)(e^v-e^{-v})}
\end{equation}
corresponding to $z_{\tau,\chi}=e^\tau\alpha^{-\frac12}$ .

Note, that when $\alpha\rightarrow 1$ we have $\chi_{top}-\chi_{bottom}\rightarrow 0$ and in the homogeneous case $\alpha=1$ there is only one turning point, as was shown to be the case in \cite{OR2}. 

As was mentioned in the introduction, unlike the homogeneous case, the point process near these turning points is not the GUE minor process. This follows immediately from the fact that we have two turning points on each edge: the interlacing property of the GUE minor process cannot hold at both of the turning points. We study the point process near these turning points in Section \ref{sec:turningPoints}.

\subsection{Triangular floor}
In the limit $u\rightarrow v$ the bounded pentagonal floor degenerates into a triangle. Differentiating \eqref{eq:spBounded} with $v=u$ we obtain
\begin{align}
\label{eq:sppTriangle}
\frac{d}{dz}\left( z\frac{d}{dz}S_{\tau,\chi}(z)\right)
=&\frac12 \frac{1}{z-e^{- u}\alpha^{-\frac 12}}
-\frac12 \frac{1}{z-e^{\tau}\alpha^{-\frac 12}}
-\frac12 \frac{1}{z-e^{\tau}\alpha^{\frac 12}}
\\\nonumber&+\frac12 \frac{1}{z-e^{ u}\alpha^{\frac 12}}.
\end{align}
The equation
\begin{equation}
\label{eq:spp=0}
\frac{d}{dz}\left( z\frac{d}{dz}S_{\tau,\chi}(z)\right)=0
\end{equation}
is again quadratic in $e^{\tau}$ and arguments similar to those used in the unbounded case give us that only one of the solutions produces a real pair $(\tau(z),\chi(z))$. Thus, as before, we obtain a parametrization of the frozen boundary. 

From \eqref{eq:sppTriangle} it is easy to see, that for fixed $u$ and $\alpha>1$ there exists a positive constant $c$ such that 
\begin{equation*}
\min_{\tau\in[-u,u]}\left\{|z_{\tau,\chi}-e^{- u}\alpha^{-\frac 12}|,
|z_{\tau,\chi}-e^{\tau}\alpha^{-\frac 12}|,
|z_{\tau,\chi}-e^{\tau}\alpha^{\frac 12}|,
|z_{\tau,\chi}-e^{ u}\alpha^{\frac 12}|\right\}>c,
\end{equation*}
which together with \eqref{eq:spBounded} imply $\chi(z)$ is bounded and the curve $(\tau(z),\chi(z))$ is a simple closed curve. This implies, in particular, that for any $\alpha>1$ the disordered region is bounded, as demonstrated in the exact sample in Figure \ref{fig:triangleSample}. Note, however, that in the homogeneous case $\alpha=1$ when $\tau(z)\rightarrow\pm u$, we have $\chi(z)\rightarrow\infty$ and the disordered region is infinite. This was first observed in \cite{BMRT}. 

As $u\rightarrow v$, we have $\chi_{top}\rightarrow \infty$, so there is only one turning point near each wall $\tau=\pm v$ when $u=v$. Moreover,
\begin{equation*}
\lim_{\alpha\rightarrow 1}\chi_{bottom}=\infty,
\end{equation*}
and when $\alpha=1$, $u=v$, we have a ``turning point at infinity'', as was observed in \cite{BMRT}.

\section{Turning points}
\label{sec:turningPoints}
We now turn to the study of the local point process at the turning points. In this section $\tau=v$, $\chi=\chi_{bottom}$ or $\chi_{top}$ and, respectively, $z_{\tau,\chi}=e^v\alpha^{\frac12}$ or $z_{\tau,\chi}=e^v\alpha^{-\frac12}$. The turning points with $\tau=-v$ are of course of the same nature, and we will not look at them.

The vertical characteristic scale at the turning point is $r^{\frac12}$, so we introduce new vertical coordinates $\tilde{h}_i$ defined by
\begin{equation*}
h_i=\lfloor\frac{\chi}r\rfloor+\frac{\tilde{h}_i}{r^\frac12}.
\end{equation*}
We also introduce new horizontal coordinates $\hat{t}_i$ which indicate the distance from the edge:
\begin{equation*}
t_i=u^\lambda_n-\hat{t}_i=
\lfloor\frac{\tau}r\rfloor-\hat{t}_i.
\end{equation*}

We deform the contours of integration in \eqref{eq:main-corr2} according to the saddle point method to pass through the critical point $z_{\tau,\chi}$. Since the $z$-contour should not cross the poles of $\Phi_{b_\lambda}(\cdot,t_1)$ the deformation is different depending on whether we have the bottom or top turning point. 

\subsubsection{Top turning point when $\Delta t\geq 0$} 
\label{sec:topTPposDelt}
\begin{figure}[ht]
\caption{\label{fig:contoursTopDeltaPos} Deformation of contours at a turning point. The shaded
region corresponds to $\Re(S_{\tau,\chi}(z)<\Re(S_{\tau,\chi}(z_{\tau,\chi})$. The solid red (inner)
and blue(outer) contours are the original contours. The dotted red
and blue contours are the deformed contours.}
\includegraphics[width=6cm]{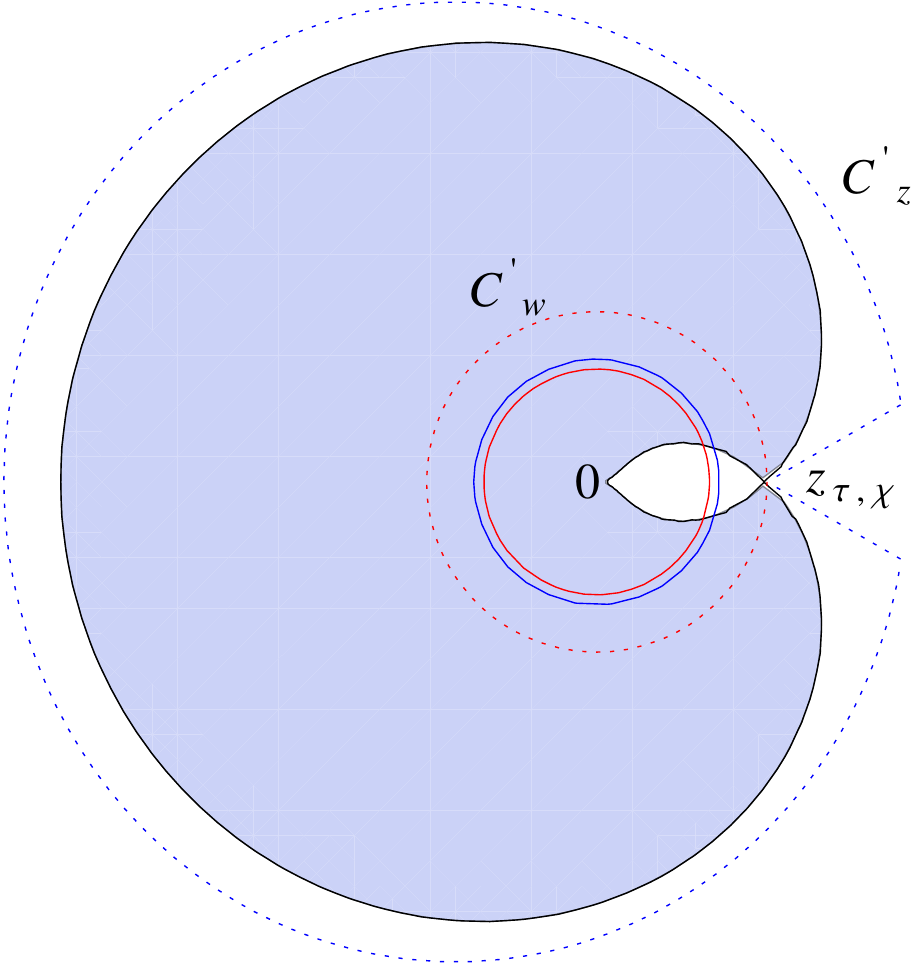}
\end{figure}

For the top turning point when $\Delta t\geq 0$ the curve $\Re S_{\tau,\chi}(z)=\Re S_{\tau,\chi}(z_{\tau,\chi})$, the original contours $C_z,C_w$ and the new contours $C'_z,C'_w$ are shown in Figure \ref{fig:contoursTopDeltaPos}. In this case the contours do not cross each other, so we do not pick up any residues as a result of the contour deformation. 

\begin{figure}[ht]
\caption{\label{fig:contourDeformedLocal} The $\zeta$ and $\omega$ contours of integration.}
\includegraphics[width=2.cm]{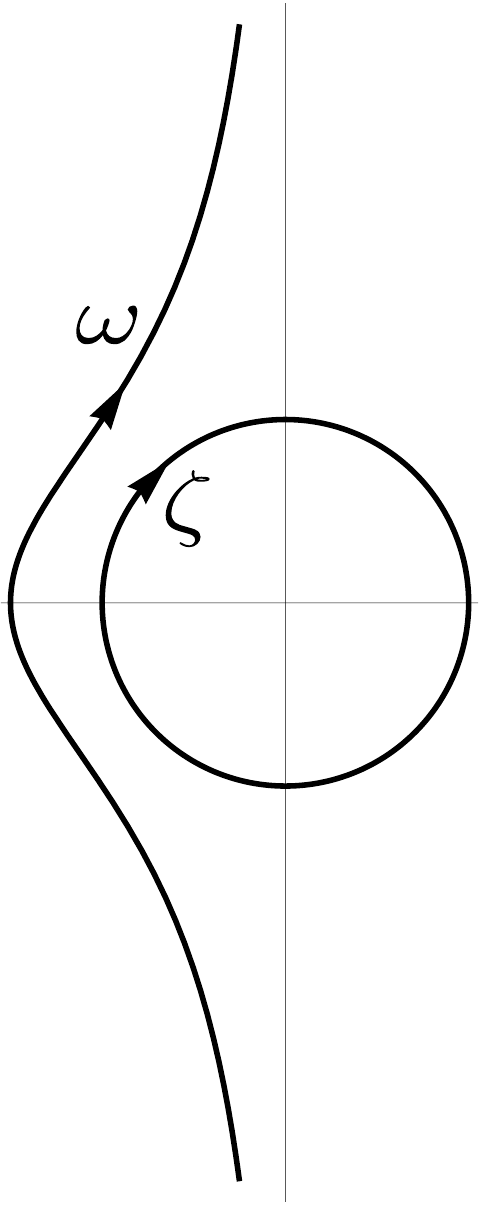}
\end{figure}

The main contribution to the integral comes from the vicinity of the critical point. 
At the turning points we have
\begin{equation*}
\frac{d}{dz}\left( z\frac{d}{dz} S_{\tau,\chi}(z)\right)\Big|_{z_{\tau,\chi}}=Q_2(z_{\tau,\chi}),
\end{equation*}
from which it is easy to see that $S''_{\tau,\chi}(z_{\tau,\chi})<0$. Changing variables of integration to $\zeta$, $\omega$ defined by
\begin{equation*}
z=z_{\tau,\chi}e^{r^{\frac12}\zeta},
\quad w=z_{\tau,\chi}e^{r^{\frac12}\omega},
\end{equation*}
we obtain
\begin{align*}
K_{\lambda_r,\bar{q}_r}((t_1,h_1),(t_2,h_2))=
&e^{\ln(z_{\tau,\chi})
\left(
  \frac{\tilde{h}_2-\tilde{h}_1}{r^{1/2}}
  -\frac{\hat{t}_2-\hat{t}_1}2
  \right)}
(1-\alpha^{-1})^{\lfloor\frac{\hat{t}_2}{2}\rfloor-\lfloor\frac{\hat{t}_1}{2}\rfloor}
\\&\times
(-r^{1/2})^{\lfloor\frac{\hat{t}_2+1}{2}\rfloor-\lfloor\frac{\hat{t}_1+1}{2}\rfloor}
\frac{r^{1/2}}{(2\pi \mathfrak{i})^2}
\\&\times
\iint e^{\frac{S''_{\tau,\chi}(z_{\tau,\chi})}2(\zeta^2-\omega^2)}
\frac{e^{\tilde{h}_2\omega}}{e^{\tilde{h}_1\zeta}}
\frac{\omega^{\lfloor\frac{\hat{t}_2+1}{2}\rfloor}}{\zeta^{\lfloor\frac{\hat{t}_1+1}{2}\rfloor}}
(1+O(r^{1/2}))\frac{d\zeta\ d\omega}{\zeta-\omega},
\end{align*}
where the contours of integration are as in Figure \ref{fig:contourDeformedLocal}.

When taking determinants of the form
\begin{equation*}
\det\left(K_{\lambda_r,\bar{q}_r}((t_i,h_i),(t_j,h_j))\right)_{1\leq i,j\leq k}
\end{equation*}
the terms
\begin{equation}
\label{eq:gaugeTurningPt}
e^{\ln(z_{\tau,\chi})
\left(
  \frac{\tilde{h}_2-\tilde{h}_1}{r^{1/2}}
  -\frac{\hat{t}_2-\hat{t}_1}2
\right)}
(1-\alpha^{-1})^{\lfloor\frac{\hat{t}_2+2}{2}\rfloor-\lfloor\frac{\hat{t}_1+2}{2}\rfloor}
(-r^{1/2})^{\lfloor\frac{\hat{t}_2+1}{2}\rfloor-\lfloor\frac{\hat{t}_1+1}{2}\rfloor}
\end{equation}
cancel out and we obtain that in the limit $r\rightarrow 0$ the leading asymptotic of the correlation kernel near the top turning point is 
\begin{equation}
\label{eq:corrKerTopTP}
\frac{r^{\frac12}}{(2\pi \mathfrak{i})^2}
\iint e^{\frac{S''_{\tau,\chi}(z_{\tau,\chi})}2(\zeta^2-\omega^2)}
\frac{e^{\tilde{h}_2\omega}}{e^{\tilde{h}_1\zeta}}
\frac{\omega^{\lfloor\frac{\hat{t}_2+1}{2}\rfloor}}{\zeta^{\lfloor\frac{\hat{t}_1+1}{2}\rfloor}}
\frac{d\zeta\ d\omega}{\zeta-\omega},
\end{equation}
where again contours of integration are as in Figure \ref{fig:contourDeformedLocal}.

\subsubsection{Bottom turning point with $\Delta t\geq 0$}
In the case of the bottom turning point, the deformed contours pass through the critical point $z_{\tau,\chi}=e^v\alpha^{\frac12}$. Since the $z$-contour should not cross the poles of $\Phi_{b_{\lambda_r},\bar{q}_r}(z,t_1)$, which are near $e^v\alpha^{-\frac12}$ or larger than $e^v\alpha^{\frac12}$, the deformed $z$-contour splits into the union of two pieces, $C'_z$ described above and $C''_z$ which is a simple closed clockwise loop around $e^v\alpha^{-\frac12}$. Moreover the $w$-contour passes over $C''_z$ during the deformation, so we pick up residues from the term $\frac1{z-w}$. To summarize, we have
\begin{multline*}
K_{\lambda,\bar{q}}((t_1,h_1),(t_2,h_2))
=\frac{1}{(2\pi \mathfrak{i})^2}
\int_{z\in C'_z}\int_{w\in C'_w}\dots+
\frac{1}{(2\pi \mathfrak{i})^2}\int_{z\in C''_z}\int_{w\in C'_w}\dots
\\+
\frac{1}{2\pi \mathfrak{i}}
\int_{z\in C''_z}
\prod_{\stackrel{t_2<m<t_1}{m+\frac12\in2\mathbb{Z}}}(1-z^{-1}x^-_m)
\prod_{\stackrel{t_2<m<t_1}{m-\frac12\in2\mathbb{Z}}}(1-zx^+_m)
z^{h_2-h_1+\frac 12 (b_\lambda(t_1)-b_\lambda(t_2))-1}dz,
\end{multline*}
where dots stand for the same integrand as in \eqref{eq:main-corr2}. The second term is exponentially small when $r\rightarrow 0$ since along the contours $C''_z$ and $C'_w$ we have
\begin{equation*}
\Re S_{\tau,\chi}(z)< \Re S_{\tau,\chi}(w),
\end{equation*}
the last term is zero since the integrand is regular at $e^v\alpha^{-\frac12}$, and the first term can be analysed as in Section \ref{sec:topTPposDelt}, giving the same result as \eqref{eq:corrKerTopTP} with $\lfloor\frac{\hat{t}_i+1}{2}\rfloor$ replaced with $\lfloor\frac{\hat{t}_i+2}{2}\rfloor$.

\subsubsection{Either turning point with $\Delta t<0$} When $\Delta t<0$, during the deformation of contours the $z$-contour passes over the $w$-contour, and we pick up residues from the term $\frac1{z-w}$ along a simple closed clockwise curve around $0$. The residues are equal to
\begin{equation*}
\frac{c}{2\pi\mathfrak{i}}
\int_{\mathfrak{i}\infty}^{-\mathfrak{i}\infty} \frac{e^{(\tilde{h}_2-\tilde{h}_1)\zeta}}
{\zeta^{\lfloor\frac{\hat{t}_1+e}{2}\rfloor-\lfloor\frac{\hat{t}_1+e}{2}\rfloor}}d\zeta,
\end{equation*}
where the contour of integration crosses the real line in the interval $(-\infty,0)$, $c$ is the gauge term \eqref{eq:gaugeTurningPt} and $e$ is $2$ for the bottom turning point and $1$ for the top one. Following \cite{OR3} we can interpret the residue term as the difference between two different expansions of
\begin{equation*}
\frac{\omega^{\lfloor\frac{\hat{t}_2+e}{2}\rfloor}}{\zeta^{\lfloor\frac{\hat{t}_1+e}{2}\rfloor}}
\frac{1}{\zeta-\omega}
\end{equation*}
depending on the sign of $\lfloor\frac{\hat{t}_1+e}{2}\rfloor-\lfloor\frac{\hat{t}_2+e}{2}\rfloor$, and thus incorporate the residue term into the main term. 

\subsubsection{The local point process at the turning points}
Combining the above results, we obtain the following theorem.

\begin{theorem}
\label{thm:corkerTurningPoint}
Let $(\tau, \chi)$ be a turning point with $\tau=\pm v$ and $\chi$ given by \eqref{eq:chiBottom} or \eqref{eq:chiTop}. Let
\begin{equation*}
t_i=\lfloor\frac{\tau}r\rfloor-\hat{t}_i,
\end{equation*}
and 
\begin{equation*}
h_i=\lfloor\frac{\chi}r\rfloor+\frac{\tilde{h}_i}{r^\frac12}.
\end{equation*}
If $\lfloor\frac{\tau}r\rfloor$ is odd, then the correlation functions near a turning point $(\tau, \chi)$ of the system with periodic weights \eqref{eq:periodicweights} are given by
\begin{equation*}
\lim_{r\rightarrow 0}r^{-\frac12} K_{\lambda,\bar{q}}((t_1,h_1),(t_2,h_2))
=\frac{1}{(2\pi \mathfrak{i})^2}
\iint e^{\frac{S''_{\tau,\chi}(z_{\tau,\chi})}2(\zeta^2-\omega^2)}
\frac{e^{\tilde{h}_2\omega}}{e^{\tilde{h}_1\zeta}}
\frac{\omega^{\lfloor\frac{\hat{t}_2+e}{2}\rfloor}}{\zeta^{\lfloor\frac{\hat{t}_1+e}{2}\rfloor}}
\frac{d\zeta\ d\omega}{\zeta-\omega},
\end{equation*}
where the contours of integration are as in Figure \ref{fig:contourDeformedLocal} and $e$ is $1$ when $\chi=\chi_{top}$ and $2$ when $\chi=\chi_{bottom}$. When $\lfloor\frac{\tau}r\rfloor$ is even, $e$ is replaced by $2-e$.
\end{theorem}
\begin{remark}
If we restrict the process to horizontal lozenges of only  even or only odd distances from the edge, then the correlation kernel in Theorem \ref{thm:corkerTurningPoint} coincides with the correlation kernel at a homogeneous turning point, obtained in \cite{OR3}. In particular, the point process of horizontal lozenges restricted to a distance of fixed parity from the edge converges to the GUE minor process.
\end{remark}

\section{Intermediate regime}
\label{sec:intermediateWeights}
In this section we study turning points in the scaling limit of random skew plane partitions under the measure \eqref{eq:qinhomog} when the weights $q_t$ are given by \eqref{eq:intermperiodicweights}. Since the nature of turning points does not depend on boundary conditions, we consider the simplest boundary which gives rise to turning points. Namely, we take $\lambda_r$ to be a staircase as before, but so that the horizontal section grows at a scale between $\frac1r$ and $\frac1{\sqrt{r}}$ and in the scaling limit the back wall converges to $V(\tau)=-|\tau|$, with $-v\leq\tau\leq v$. 

The function $S_{\tau,\chi}$ is the same as in the case of homogeneous weights and is given by \eqref{eq:S} with $\alpha=1$. It follows that the macroscopic behaviour of the system is the same as in the homogeneous case. In particular it has the same frozen boundary, and only one turning point near each vertical $\tau=\pm v$. Moreover, the correlation kernel in the bulk matches with the homogeneous case and in the bulk horizontal lozenges are distributed according to a $\mathbb{Z}\times\mathbb{Z}$ translation invariant ergodic Gibbs measure. In particular, periodicity disappears in the scaling limit in the bulk. However, even tough there is only one turning point near each of the walls $\tau=\pm v$, the local point process is not the GUE minor process. Following the method from Section \ref{sec:turningPoints} we obtain the following theorem.

\begin{theorem}
\label{thm:corkerTurningPointIntermediate}
Let $(\tau, \chi)$ be a turning point with $\tau=\pm v$. Let
\begin{equation*}
t_i=\lfloor\frac{\tau}r\rfloor-\hat{t}_i,
\end{equation*}
and 
\begin{equation*}
h_i=\lfloor\frac{\chi}r\rfloor+\frac{\tilde{h}_i}{r^\frac12}.
\end{equation*}
If $\lfloor\frac{\tau}r\rfloor$ is odd, then the correlation functions near a turning point $(\tau, \chi)$ of the system with periodic weights \eqref{eq:intermperiodicweights} are given by
\begin{multline*}
\lim_{r\rightarrow 0}r^{-\frac12} K_{\lambda,\bar{q}}((t_1,h_1),(t_2,h_2))
\\=\frac{1}{(2\pi \mathfrak{i})^2}
\iint e^{\frac{S''_{\tau,\chi}(z_{\tau,\chi})}2(\zeta^2-\omega^2)}
\frac{e^{\tilde{h}_2\omega}}{e^{\tilde{h}_1\zeta}}
\frac{\omega^{\lfloor\frac{\hat{t}_2+1}{2}\rfloor}}{\zeta^{\lfloor\frac{\hat{t}_1+1}{2}\rfloor}}
\frac{(\omega-\gamma)^{\lfloor\frac{\hat{t}_2+2}{2}\rfloor}}{(\zeta-\gamma)^{\lfloor\frac{\hat{t}_1+2}{2}\rfloor}}
\frac{d\zeta\ d\omega}{\zeta-\omega},
\end{multline*}
where the contours of integration are as in Figure \ref{fig:contourDeformedLocal} with the $\zeta$-contour containing both $0$ and $\gamma$. If $\lfloor\frac{\tau}r\rfloor$ is even, then the exponents of $\zeta$ and $\omega$ switch places with those of $\zeta-\gamma$ and $\omega-\gamma$.
\end{theorem}

\bibliography{mybib}
\bibliographystyle{alpha}

\end{document}